\documentclass[12pt]{extarticle}
\usepackage{amsmath, amsthm, amssymb, mathtools, hyperref,color}
\usepackage{graphicx}

\newtheorem{theorem}{Theorem}
\numberwithin{theorem}{section}
\newtheorem{proposition}[theorem]{Proposition}
\newtheorem{proposition-definition}[theorem]{Proposition-Definition}
\newtheorem{lemma}[theorem]{Lemma}
\newtheorem{corollary}[theorem]{Corollary}
\newtheorem{definition}[theorem]{Definition}
\newtheorem{remark}[theorem]{Remark}

\newtheorem{example}[theorem]{Example}

\newtheorem{prop-def}[theorem]{Proposition-Definition}
\newtheorem{principle}{Principle}

\title{Equations over Polyhedral Semirings}
\author{Madhusudan Manjunath\footnote{MM was supported by a MATRICS grant (MTR/2019/000462) of the Department of Science and Technology (DST), India. Part of this work was carried out while visiting ICTS, Bangalore and IHES, Bures-sur-Yvette. We thank the institutes for their generous support and warm hospitality.}}
\begin{document} 
\maketitle
\begin{abstract}
We study the theory of equations in one variable over polyhedral semirings. The article revolves around a notion of solution to a polynomial equation over a polyhedral semiring. 
Our main results are a characterisation of local solutions in terms of the coefficients, a local-global principle, and the basics of multiplicity and discriminants. Our main sources of motivation are tropical geometry and the theory of exceptional points in non-Hermitian physics. 
\end{abstract}
\section{Introduction}

Solutions to polynomials is one of the most fundamental problems in Mathematics. A closely related problem is the study of dependence of the solutions on (some aspect of) the coefficients of the underlying polynomial. One version of this problem is the following:
\begin{center}
Let $R=\mathbb{K}[x_1,\dots,x_n]$ where $\mathbb{K}$ is a field. Given a polynomial $p \in R[t]$, determine the dependance of the solutions to $p$ on the Newton polytopes of its coefficients.
\end{center}

Solutions to $p$ are described in terms of certain ``generalised" Puiseux series and hence, we may ask: How does polyhedral data, akin to the Newton polytope, associated to solutions to $p$ depend on the Newton polytopes of its coefficients?  The current article investigates this problem in terms of certain semirings associated to polyhedra referred to as \emph{polyhedral semirings}. 
The central theme is the notion of solution to polynomials over polyhedral semirings.  This article revolves around developing theory underlying this notion. This prospect of developing algebraic geometry over polyhedral semirings has been raised in this article by Speyer and Sturmfels \cite[Page 165]{SpeStu09}. This continues a theme from \cite{Man17} where the notion of a syzygy over a polytope semiring was studied.

A broad motivation for this investigation is to study fundamental themes in algebraic geometry such as syzygies and rationality in terms of their polyhedral counterparts. Another perspective is as an extension of tropical geometry, we shall discuss this in more detail later (Subsection \ref{motrel_intro}).

A polyhedral semiring $\mathcal{A}$ consists of certain polyhedra in $\mathbb{R}^n$ with convex hull $\oplus$ as the addition and Minkowski sum $\odot$ as the multiplication along with an additive identity $0_{\mathcal{A}}$. A polyhedral polynomial $\phi$ in one variable with coefficients in $\mathcal{A}$ is a symbol of the form $Q_d \odot Y^d \oplus \dots \oplus Q_0$ where each $Q_i \in \mathcal{A}$. A polyhedron $P$ is called a \emph{root} of $\phi$ (or a solution to $\phi$) if every vertex of the evaluation $\phi(P)$ of $\phi$ at $P$ is shared by at least two summands.  This definition is inspired by its counterpart in tropical geometry and the notion of a syzygy over a polytope semiring.  We refer to Section \ref{polyhed_sect} for details.   We summarise our main results in the following.

\begin{itemize} 
\item {\bf Local Solutions}: This part deals with solutions that are ``local'' in the sense that their normal cone is contained in the normal cone of a fixed vertex of the Minkowski sum of the coefficients of $\phi$. We characterise, under a genericity assumption on $\phi$, local polyhedral solutions to $\phi$ in terms of a certain (labelled) fan constructed from its coefficients.

\item{\bf Local-Global Principle}: This part involves gluing together local solutions to obtain a global one (a solution with ``maximum'' normal fan). 
Inspired by local-global principles in number theory,  we characterise collections of local solutions that can glued to a global solution.

\item{\bf Discriminants}: This section develops the notion of multiplicity and discriminant in the polyhedral setting. Our main result here is a polyhedral analog of the characterisation of degenerate  polynomials in terms of the vanishing of their coefficients on the corresponding discriminant. 

\end{itemize}

We elaborate each topic to state the corresponding theorems.

\subsection{Local Solutions}

 Unlike the classical case where every non-zero polynomial in one variable has finite many solutions, all polyhedral polynomials have infinitely many solutions.  In particular,  if $P$ is a solution to $\phi$, then so is $P \odot C$ for any polyhedral cone $C$. Keeping this is mind, we focus on certain solutions that we refer to as \emph{Minkowski-Weyl minimal} solutions: a solution $P$ is called Minkowski-Weyl minimal with respect its vertex set if its recession cone is minimal among all solutions with the same vertex set as $P$.  
 The existence of Minkowski-Weyl minimal solutions is not evident a priori.

A starting observation is that any vertex of any root of $\phi$ is of the form $\rho^{(\nu)}_{i,j}:=-(\nu_i-\nu_j)/(i-j)$ where $i \neq j$, $\nu$ is a vertex of $M$ and $\nu=\nu_{i_0}+\nu_{i_1}+\dots+\nu_{i_{r+1}}$ is its Minkowski decomposition (Proposition \ref{mindecgenpol_prop}). This raises the problem of picking subsets of vertices of the form $\rho^{(\nu)}_{i,j}$ and associating to each element an appropriate cone that serves as its ``normal cone" so that this data together yields a solution to $\phi$.  This problem seems rather involved at this level of generality, as we shall see further. Hence, we restrict to \emph{local} solutions.

Fix a vertex $v$ of $M$, a $v$-local solution to $\phi$ is a solution whose normal fan is contained in the normal cone of $v$.  A $v$-local solution $P$ is called Minkowski-Weyl minimal if its recession cone is minimal among all $v$-local solutions with the same vertex set.  Our main theorem characterises Minkowski-Weyl minimal $v$-local solutions to $\phi$ in terms of objects called \emph{local compatible systems (LCS)} constructed in terms of the coefficients of $\phi$.  In the following, we provide an informal description of LCS. 

We start by refining the normal cone of $v$ to a fan $\mathcal{F}_v$ as follows. Roughly speaking, the maximal cones of $\mathcal{F}_v$ correspond to candidate normal cones of vertices of a  $v$-local solution. We label the maximal cones of $\mathcal{F}_v$ with two types of data, namely \emph{the primary labels} and \emph{the secondary labels}.  For a maximal cone $C$ of $\mathcal{F}_v$,  its primary label keeps track of points $\rho^{(v)}_{i,j}$ such that the affine cone $\rho^{(v)}_{i,j}+C^{\star}$, where $C^{\star}$ is the dual of $C$, is a solution to $\phi$.  On the other hand, the secondary labels of $C$ keep track of pairs of indices $((i,j),(\tilde{i},\tilde{j}))$ such that $\rho^{(v)}_{i,j}$ and $\rho^{(v)}_{\tilde{i},\tilde{j}}$ are candidate vertices (among possibly others) of a $v$-local solution with $\rho^{(v)}_{i,j}$ as a vertex whose normal cone is $C$ .  A local compatible system at $v$ consists of data of the form $(\mathcal{C},\mathcal{I})$ where  $\mathcal{C}=\{C_k\}_{k=1}^{t}$ is a collection of maximal cones of $\mathcal{F}_v$ and $\mathcal{I}=\{(i_k,j_k)\}_{k=1}^{t}$ is a corresponding collection of pairs of integers that together satisfy certain compatibility conditions (Definition \ref{lcs_def}).  A local compatible system ``glues'' maximal cones of $\mathcal{F}_v$ to yield $v$-local solutions to $\phi$. Naively, one may expect local compatible systems to be in bijection with  \emph{Minkowski-Weyl minimal}  polyhedral solutions to $\phi$. But this is not quite the case, there is an additional ``convexity property" (namely the convexity of $\cup_{k=1}^{t}C_k$) that local compatible systems do not capture.  Taking cue from this, we consider certain generalisations of polyhedra referred to as \emph{vertex-cone collections} (Definition \ref{vercon-def}). We extend the notion of polyhedral solutions and their Minkowski-Weyl minimality to vertex-cone collections.  Under a genericity assumption on $\phi$, local compatible systems of $\phi$ are in bijection with its Minkowski-Weyl minimal vertex-cone collection solutions (Theorem \ref{genminwey_theo}). We then take the additional convexity property into account and associate polyhedra to certain polyhedral convex subsets of $\cup_{k=1}^{t}C_k$  to show the following.

\begin{theorem}\label{minweychar_theointro} Let $\phi_{\rm gen}$ be a generic polyhedral polynomial. Fix a vertex $v$ of $M$. A polyhedron $P_0$ with vertex set $V$  is a $v$-local Minkowski-Weyl minimal solution to $\phi_{\rm gen}$ if and only if there is a local compatible system $(\mathcal{C},\mathcal{I})$ at $v$ where $\mathcal{C}=\{C_k\}_{k=1}^{t}$ and $\mathcal{I}=\{(i_k,j_k)\}_{k=1}^{t}$ such that $V=\{\rho^{(v)}_{i_k,j_k}\}_{k=1}^t$ and  $P_0$ is the associated polyhedron of the restriction of $(\mathcal{C},\mathcal{I})$ to a maximal polyhedral convex subset of $\cup_{k=1}^{t}C_k$. 
 \end{theorem}
 
 {\bf Main Ingredients of the Proof:} Construction of $v$-local polyhedral solutions involves two phenomena. The first phenomenon involves gluing maximal cones of the fan $\mathcal{F}_v$ by keeping track of primary and secondary labels across these maximal cones.  This is reminiscient of wall crossing phenomena in algebraic geometry \cite{BayMac22}. The second phenomenon involves transforming the result of the gluing construction so that it is the normal fan of a polyhedron. This involves ``extracting'' maximal polyhedral convex subsets from the result of the first step. The proof of Theorem \ref{genminwey_theo} implements the first step.  Taking Theorem \ref{genminwey_theo} as input, the proof of Theorem \ref{minweychar_theointro} executes the second step.

 The main technical difficulties are: i. to show that any $v$-local  Minkowski-Weyl minimal vertex-cone collection solution corresponds to a local compatible system and ii. to associate a local compatible system  to any ($v$-local) Minkowski-Weyl minimal polyhedral solution. The first difficulty is handled via the genericity assumption and two continuation principles (Principles \ref{firstcont_princ} and \ref{seccont_princ}). The second difficulty is handled via a notion of completion that associates a Minkowski-Weyl minimal vertex-cone collection solution
 to any Minkowski-Weyl minimal polyhedral solution, both $v$-local (Subsection \ref{comasspshr_subsect}).

 {\bf The Non-Local Case:} The local hypothesis along with the genericity assumption on the polyhedral polynomial allow the characterisation of Minkowski-Weyl minimal vertex-cone collection solutions. Without either of these assumptions, such a characterisation seems very involved as the following example indicates. 
 
 Consider a three-polytope $Q_1$ with adjacent vertices $v_1,v_2$ such that the union of their normal cones is not convex. For example, any pair of adjacent vertices of a regular icosahedron has  this property. Let $Q_0$ be the edge formed by $v_1$ and $v_2$. Consider the polyhedral polynomial $Q_1 \odot Y \oplus Q_0$. The polytope $M=Q_1 \odot Q_0$ has vertices $2 v_1$ and $2v_2$, and they form an edge. The points $\rho^{(2 v_1)}_{0,1}=\rho^{(2 v_2)}_{0,1}=0$. Normal fans of non-local solutions to this polynomial can contain convex subsets of the union of the normal cones of $2v_1$ and $2v_2$ glued with possibly other cones.  Their characterisation seems difficult. We refer to Section \ref{example_sect} for difficulties in the local case when the genericity assumption is dropped.

\subsection{Local-Global Principle}
 
 Generally speaking, local-global principles describe ``global'' phenomena in terms of corresponding ``local'' phenomena where ``global'' and ``local'' are defined according to the context. The most basic version of this principle states that a rational quadratic form has a non-trivial rational solution if and only if it has a non-trivial solution over $\mathbb{Q}_p$ for every prime $p$ and a non-trivial solution over $\mathbb{R}$. An analogous principle does not hold for forms of higher degree, starting from the cubics and there are further obstructions to the existence of global solutions. We refer to \cite{Par14} for a survey.  Inspired by this, we seek a ``local-global principle'' in this polyhedral setting.  We formulate the following question.

\begin{center}
Suppose that a polyhedral polynomial has a complete $v$-local solution, i.e. with normal cone $N_M(v)$, for each $v \in M$. Under what additional conditions, if any, does it have a global solution, i.e. with normal cone containing that of $M$? 
 \end{center}
 
 We identify a simple ``orientation'' condition and show that it suffices to glue local solutions to a global one. 
 
 \begin{theorem}{\rm {(\bf Local-Global Principle})}\label{locglo_theo}  A polyhedral polynomial $\phi$  has a global solution if and only if 
 for each vertex $v$, there is a complete $v$-local solution $S_v$ such that the following property holds: 
\begin{enumerate}
\item[(P)] \label{ori_prop} If $N_M(v) \cap N_{P_0}(\gamma)$ is full-dimensional, then $\gamma$ is a vertex of $S_v$
\end{enumerate}

where $P_0:={\rm CH}(\cup_{v \in V(M)}V(S_v))$ (here, ${\rm CH}(.)$ is the convex hull) and  $\gamma$ is a vertex of $P_0$.  Furthermore, the polytope $P_0$ is a global solution to $\phi$. 
 \end{theorem}

If $\phi$ is linear, i.e. $\phi=Q_1 \odot Y \oplus Q_0$ and if $Q_1,~Q_0$ are polytopes, then it has a global solution if and only if  $Q_1$ is a Minkowski summand of $Q_0$. 
In this case,  Theorem \ref{locglo_theo} is a characterisation of Minkowski summands. In Section \ref{locglo_sect}, we use this connection to reprove a criterion for weak Minkowski summands due to Shephard. An interesting future direction is to view Theorem \ref{locglo_theo} as a ``shadow'' of its counterpart over algebraic varieties \cite{Par14}. 

 
\subsection{Multiplicity and Discriminants}  
 
 A satisfactory theory of multiplicity and discriminants is central to any theory of polynomial equations. This part draws inspiration from both the classical and tropical theory of discriminants. 
 We define the notion of multiplicity of (polyhedral) roots and via this, a notion of degeneracy and its weak counterpart for polyhedral polynomials. 
 This leads us to the following definition of polyhedral discriminants. Let $\Xi$ be a finite subset of $\mathbb{Z}_{\geq 0}$, a polyhedral polynomial $\Delta_{\Xi} \in \mathcal{A}[Y_i|~i \in \Xi]$   is called a \emph{polyhedral discriminant} with respect to $\Xi$ if the following holds: $\Delta_{\Xi}$ vanishes on a tuple if and only if the polynomial in $\mathcal{A}[Y]$  with this tuple of coefficients is degenerate and is in a certain form called \emph{product form} (Definition \ref{polydisc_def}). 
 Via a transformation called \emph{polyhedralisation} that takes polynomials to polyhedral polynomials (Section \ref{polyhedralise_sect}), we consider \emph{polyhedralised discriminants} that are polyhedralisations of (classical) discriminants.  Our main result is the following.
 
 \begin{theorem}\label{polydisc_theo}
 Let $\phi=\oplus_{i \in \Xi} (Q_i \odot Y^i) \in \mathcal{A}[Y]$ be in product form.  If $\phi$ is degenerate, then $(Q_i|~i \in \Xi)$  is a root of the polyhedralised discriminant $\tilde{\Delta}_{\Xi}$. Conversely, if the tuple $(Q_i|~i \in \Xi)$ is a root of $\tilde{\Delta}_{\Xi}$,  then $\phi$  is weakly degenerate.  
 \end{theorem}
 
 {\bf Main Steps of the Proof:} The proof builds on the fact that the analogous statement holds in the tropical setting.  The polyhedral polynomial $\phi$ localises to a collection of univariate tropical polynomials each in product form and hence, the tuple of coefficients of each of these tropical polynomials is a root of the corresponding tropicalised discriminant.  This implies that $(Q_i|~i \in \Xi)$  is a root of $\tilde{\Delta}_{\Xi}$. We refer to Proposition \ref{polydisc_prop} for details.    Conversely, the polyhedralised discriminant ``localises" to a collection of tropicalised discriminants and the resulting tropical solutions of high multiplicity are lifted to a polyhedral one (Proposition \ref{discpoly_prop}).

An attempt to strengthen the second part of Theorem \ref{polydisc_theo} to $\phi$ being degenerate leads to interesting challenges in gluing cones and is a topic for further work.




\subsection{Motivation and Related Work}\label{motrel_intro}

This investigation is inspired by the following sources:

\begin{enumerate}

\item {\bf Ramification of Hypersurfaces and Exceptional Points}:  The study of ramification of hypersurfaces under a projection map onto a coordinate plane arises in the study of exceptional points in non-Hermitian physics \cite{Hei12}.  A typical situation is as follows: the projection map is ramified at a point (the origin, say) and the set of lines through the origin under whose restriction the projection map remains ramified  is important. In a recent work \cite{BJMN23}, tropical geometric methods were used to study this problem in the planar case. The current article was partially motivated by generalisation of this tropical approach to higher dimensions and is a topic for future work.



\item {\bf Higher-Rank Tropical Geometry}:  Already in the initial conception of tropical geometry \cite{SpeStu09}, Speyer and Sturmfels foresaw  the possibility of ``algebraic geometry" over semirings in the same flavour as the tropical semiring such as the polyhedral semiring and suggested this as a research direction.  Currently,  tropical geometry over fields with higher rank valuations is taking shape, \cite{FosRan16, JosSmi23,AmiIri22}.
The definition of higher rank valuations underlies a well-order on the value group, for instance the lexicographic ordering on Euclidean space.  
Instead of picking a well-order, if we deal with all ``possible term orders'' at once by considering polyhedral approximations of ${\rm New}(s)$ where $s$ is the underlying power series, then we land in the setting of polyhedral semirings. The current article fits into this framework. This also suggests studying semiring homomorphisms from the function field of a variety to a polyhedral semiring.  
 
\end{enumerate}




{\bf Acknowledgements:}  We thank Ayan Banerjee, Rimika Jaiswal and Awadhesh Narayan for their collaboration on exceptional points. We thank Lorenzo Fantini for interesting discussions on polyhedral semirings while the author was visiting IHES.

\section{Polyhedral Semirings}\label{polyhed_sect}

The main content of this section is the definition of a root of a polyhedral polynomial. 
We first define certain basic semirings consisting of polyhedra that we refer to as \emph{$\omega$-positive semirings} and define an arbitrary polyhedral semiring as a subring of such an $\omega$-positive semiring.   Fix a vector $\omega$ in $(\mathbb{R}^n)^{\star}$ whose coordinates are rationally independent, i.e. $\mathbb{Q}$-linearly independent.   

A polyhedron $P$ in $\mathbb{R}^n$ is called $\omega$-positive if its inner normal fan contains the vector $\omega$. Equivalently, $P$ is $\omega$-positive if the linear program  ${\rm min}_{{\bf x} \in P}~\omega({\bf x})$ is bounded.  We denote the Minkowski-Weyl cone of $P$, also known as the recession cone, by ${\rm MWC}(P)$.  




\begin{proposition-definition}{\rm {(\bf $\omega$-Positive Semirings})}
 The set  $\mathcal{P}_\omega$ of all $\omega$-positive rational polyhedra in $\mathbb{R}^n$ is closed under the operations convex hull $\oplus$ and Minkowski sum $\odot$. This set along with a distinguished element $0_{\mathcal{P}_\omega}$ as the additive identity forms a semiring with $\oplus$ and $\odot$ as the addition and multiplication, respectively. 
\end{proposition-definition}

\begin{remark}
\rm Note that $P \odot 0_{\mathcal{P}_\omega}=0_{\mathcal{P}_\omega}$ for all $P \in \mathcal{P}_\omega$. To see this, note that $P \odot (Q \oplus 0_{\mathcal{P}_\omega})=P \odot Q$ for all $Q \in \mathcal{P}_\omega$ and choosing $Q$ to be points over $\mathbb{R}^n$ yields the claim. \qed
\end{remark}


\begin{remark}
\rm Our notion of $\omega$-positive semirings is inspired by the work of Aroca and Ilardi \cite{AroIla09} where they introduce the notion of $\omega$-positive (multivariate) Puiseux series and show that for a fixed $\omega \in (\mathbb{R}^n)^{\star}$ with rationally independent coordinates, the set of all $\omega$-positive Puiseux series (in $n$ variables) is algebraically closed. According to them, an $n$-variate Puiseux series is called $\omega$-positive if there is a translate of the set of exponents (of monomials in its support) that is contained in a strongly convex, rational polyhedral cone and this cone is, in turn, contained in the half-space $\{{\bf x}|~\omega({\bf x}) \geq 0)\}$. It follows that a polyhedron is the Newton polyhedron of an $\omega$-positive Puiseux series if and only if it is $\omega$-positive. We also refer to \cite[Section 3]{Val22} for a treatment of polyhedral semirings.
 \qed
\end{remark}


We refer to the book of Golan \cite{Gol99} for a theory of semirings, in general.  The Bachelor's thesis \cite{Val22} is also a handy reference for polyhedral semirings.

\begin{definition}{\rm {(\bf Polyhedral Semirings})}
A semiring is called a polyhedral semiring if it contains an additive identity, a multiplicative identity and admits a semiring endomorphism to an $\omega$-positive semiring.
\end{definition}



Let $\mathcal{A}$ be a polyhedral semiring with additive identity $0_{\mathcal{A}}$ and embedded in $\mathcal{P}_\omega$ for a suitable $\omega$.  A non-zero polynomial in the indeterminates $Y_1,\dots,Y_h$ with coefficients in $\mathcal{A}$ is a symbol of the form $\oplus_{\alpha \in \mathcal{S}}(Q_\alpha \odot Y_1^{\alpha_1} \odot \dots \odot Y_h^{\alpha_h})$ where $\alpha=(\alpha_1,\dots,\alpha_h), \mathcal{S}$ is a non-empty, finite subset of $(\mathbb{Z}_{\geq 0})^h$ and $Q_{\alpha} \neq 0_{\mathcal{A}}$ is in $\mathcal{P}_\omega$  for all $\alpha \in \mathcal{S}$.  We use mult-index notation and denote $Y_1^{\alpha_1} \odot \dots \odot Y_h^{\alpha_h}$ by $Y^{\alpha}$.


The addition and multiplication operations on $\mathcal{A}$ extend to the set of all polynomials in  $Y_1,\dots,Y_h$ via the additional relation $Y^{\alpha^{(1)}} \odot Y^{\alpha^{(2)}}=Y^{\alpha^{(1)}+\alpha^{(2)}}$.

 \begin{proposition}
 The set of all polynomials in the indeterminates $Y_1,\dots,Y_h$ with coefficients in $\mathcal{A}$ is a semiring with respect to the operations $\oplus$ and $\odot$ with additive identity $0_{\mathcal{A}}$.
  \end{proposition}
 
 
  Let $\phi= \oplus_{\alpha \in \mathcal{S}}(Q_\alpha \odot Y_1^{\alpha_1} \odot \dots \odot Y_h^{\alpha_h}) \in \mathcal{A}[Y_1,\dots,Y_h]$ and $P=(P_1,\dots,P_h)$ be an $h$-tuple of $\omega$-positive polyhedra. The \emph {evaluation} $\phi(P)$ of $\phi$ at $P$ is the polyhedron $\oplus_{\alpha \in \mathcal{S}}(Q_\alpha \odot P_1^{\odot \alpha_1} \odot \dots \odot P_h^{ \odot \alpha_h})$  where $P_j^{\odot \alpha_j}$ is the Minkowski sum $P_j \underbrace{\odot \cdots \odot}_{\alpha_j \text{ times }} P_j$ if $\alpha_j$ is a positive integer and $P^{\odot 0}$ is the multiplicative identity of $\mathcal{A}$.  Note that the polyhedron $\phi(P)$ and every summand $Q_\alpha \odot P_1^{\odot \alpha_1} \odot \dots \odot P_h^{ \odot \alpha_h}$ is $\omega$-positive.  



\begin{definition}{\rm ({\bf Root of a Polyhedral Polynomial})}\label{polroot_eq}
An $h$-tuple of $\omega$-positive polyhedra $P_0=(P_{0,1}\dots P_{0,h})$  is said to be a \emph{solution} to $\phi$ if one of the following conditions is satisfied:
\begin{enumerate}

\item The evaluation  $\phi(P_0)$ is $0_{\mathcal{A}}$.

\item  Otherwise, every vertex of $\phi(P_0)$ is shared by at least two distinct summands, i.e. for every vertex $v$ of $\phi(P_0)$, there exist summands $Q_l \odot P_0^{\odot l}$ and $Q_m \odot P_0^{\odot m}$ of $\phi(P_0)$ with $l \neq m$ such that $v$ is a vertex of both these summands. 
 \end{enumerate}
\end{definition}

Definition \ref{polroot_eq} is justified by the following properties.
\begin{enumerate}

\item  Let $p \in \mathbb{K}[x_1,\dots,x_h]$ where $\mathbb{K}$ is the field of $\omega$-positive Puiseux series.  For any root $s$ of $p$,  the Newton polyhedron of $s$ is a root of its ``polyhedralisation'' ${\rm polyhed}(p) \in  \mathcal{P}_\omega[Y_1,\dots,Y_h]$.   We refer to Proposition \ref{rootpres_prop} for more details.

\item  This definition of a root generalises the corresponding notion over the tropical semiring.   As explained in \cite[Page 165]{SpeStu09}, the tropical semiring $\mathbb{T}=(\mathbb{R} \cup \{ \infty\}, {\rm min},+)$  embeds into the one-dimensional $\omega$-positive semiring (where $\omega$  is any positive integer) via the semiring endomorphism $\iota$ that takes  $a$ to $[a,\infty]$. Furthermore, a real number $c$ is a root of the tropical polynomial ${\rm min}_{j \in S \subseteq [0,\dots,n]}\{a_j+ j \cdot x \}$ where  $n \in \mathbb{Z}_{>0}$  and $a_j \in \mathbb{R}$  if and only if $\iota(c)$ is a root of the polynomial $\oplus_{ j \in S} (\iota(a_j) \odot Y^{\odot j})$ in $\mathcal{A}_{\mathbb{T}}[Y]$.



\end{enumerate}

\begin{remark}\label{minsum_rem}
\rm
A linear polyhedral polynomial $Q_1 \odot Y \oplus Q_0$ where $Q_1$ and $Q_0$ are polytopes has a polytopal solution if and only if $Q_1$ is a Minkowski summand of $Q_0$. The set of all  ``weak Minkowski summands'' of a polytope forms a cone called the \emph{type cone} and has been studied in literature \cite{CaDoGoRoYi22}. 
\qed
\end{remark}






\section{Polyhedralisation}\label{polyhedralise_sect}

We associate polyhedral polynomials to certain polynomials with coefficients over an $\omega$-positive Puiseux series field.


Fix any $\omega \in (\mathbb{R}^n)^\star$ whose coordinates are rationally independent. Let $\mathbb{K}$ be the field of $\omega$-positive Puiseux series in $n$ variables.  For an element $s \in \mathbb{K}^{\star}$, let ${\rm New}(s)$ be the convex hull of the exponents of its support and let ${\rm New}(0_{\mathbb{K}})=0_{\mathcal{P}_\omega}$.  Note that in general, ${\rm New}(s)$ is a convex body and not, in general, a polyhedron. However, the subset of elements $s$ in $\mathbb{K}^{\star}$ such that ${\rm New}(s)$ is an $\omega$-positive polyhedron, along with $0_{\mathbb{K}}$ is a subring of $\mathbb{K}$. We denote this subring by $\mathbb{K}_P$.  Fix a positive integer $h$ and let $\mathbb{K}_P[x_1,\dots,x_h]$ be the polynomial ring in $h$ variables with coefficients over $\mathbb{K}_P$.  
The \emph{polyhedralisation map} ${\rm polyhed}:\mathbb{K}_P[x_1,\dots,x_h] \rightarrow \mathcal{P}_{\omega}[Y_1,\dots,Y_h]$ is defined as follows. For $p= \sum_{\alpha|c_{\alpha} \neq 0} c_{\alpha} {\bf x}^{\alpha} \in \mathbb{K}_P[x_1,\dots,x_h]$,
\begin{center}
${\rm polyhed}(p)=
\begin{cases}
 0_{\mathcal{P}_{\omega}},{\rm if~} p=0,\\
 \oplus_{\alpha}({\rm New}(c_\alpha) \odot {\bf Y}^{\odot \alpha}),{\rm otherwise}. 
\end{cases}
$
\end{center}
where ${\bf Y}^{\odot \alpha}= Y_1^{\odot \alpha_1} \odot \dots \odot Y_h^{\odot \alpha_h}$ for $\alpha=(\alpha_1,\dots,\alpha_h)$.

The semiring structure on $\mathcal{P}_{\omega}$ extends to $\mathcal{P}_{\omega}[Y_1,\dots,Y_h]$ in the standard way. The map ${\rm polyhed}$ is not quite a homomorphism of semirings, nevertheless the following holds. 
\begin{proposition}
If polynomials $p_1,p_2 \in \mathbb{K}_P[x_1,\dots,x_h]$ have disjoint support, then ${\rm polyhed}(p_1+p_2)={\rm polyhed}(p_1) \oplus {\rm polyhed}(p_2)$. The map ${\rm polyhed}$ is multiplicative, i.e. ${\rm polyhed}(p_1 \cdot p_2)={\rm polyhed}(p_1) \odot {\rm polyhed}(p_2)$ for any $p_1,p_2 \in \mathbb{K}_P[x_1,\dots,x_h]$.
\end{proposition}

 The support ${\rm Supp}(\mathcal{B})$ of a convex set $\mathcal{B}$ is the set of linear functionals $\ell$ such that ${\rm min}_{{\bf x} \in \mathcal{B}} \ell({\bf x})$ exists. 
 The following proposition shows that any root $s$ of a polynomial in $\mathbb{K}_P[x]$ is ``almost polyhedral", i.e. ${\rm New}(s)$ has finitely many vertices. 
  By a vertex of  $\mathcal{B}$, we mean a point ${\bf p}$ in $\mathcal{B}$ for which there is a linear functional $\ell \in {\rm Supp}(\mathcal{B})$ such that ${\rm min}_{{\bf x} \in \mathcal{B}} \ell({\bf x})$ has ${\bf p}$ as its unique solution.

\begin{proposition}\label{finver_prop}
Let $p=\sum_{i: c_i \neq 0} c_i x^i \neq 0_{\mathbb{K}} \in \mathbb{K}_P[x]$. For any non-zero root $s$ of $p$ such that ${\rm Supp}({\rm New}(s)) \subseteq {\rm Supp}(M)$ where $M=\odot_i {\rm New}(c_i)$, the convex set ${\rm New}(s)$ has finitely many vertices. 
\end{proposition}
\begin{proof}
Consider a vertex $w$ of ${\rm polyhed}(p)(({\rm New}(s)))$. Note that $w=v_i+i \cdot \gamma$ for vertices $v_i,\gamma$ of ${\rm New}(c_i)$ (for some $i \in {\rm Supp}(p)$) and ${\rm New}(s)$, respectively.   By definition of a vertex, there is a linear functional $\ell$ that attains a unique minimum over ${\rm polyhed}(p)({\rm New}(s))$ at $w$. This linear functional $\ell$ also attains unique minima at $v_i$ and $\gamma$ over ${\rm New}(c_i)$ and ${\rm New}(s)$, respectively.  Since the coefficient of ${\bf t}^w$ in $p(s)$ is zero, there must be another summand $j \in {\rm Supp}(p)$ apart from $i$ such that $w=v_j+j \cdot \tilde{\gamma}$ for vertices $v_j, \tilde{\gamma}$ of ${\rm New}(c_j)$ and ${\rm New}(s)$, respectively. Furthermore, we claim that $\gamma=\tilde{\gamma}$. Suppose not, then $\ell(\gamma) < \ell(\tilde{\gamma})$.  In this case, $\ell(v_j+j \cdot \gamma)<\ell(v_j+j \cdot \tilde{\gamma})$. Since $v_j+j \cdot \gamma \in {\rm polyhed}(p)(({\rm New}(s)))$, this contradicts the fact that $\ell$ attains a unique minimum at $w$ over ${\rm polyhed}(p)(({\rm New}(s)))$. Hence, $w=v_i+i \cdot \gamma=v_j+j \cdot \gamma$ and $\gamma=-(v_i-v_j)/(i-j)$. 
Furthermore, every vertex $\tilde{\gamma}$ of ${\rm New}(s)$ is uniquely minimised by a linear functional $\ell \in {\rm Supp}({\rm New}(s))$. Since ${\rm Supp}({\rm New}(s))={\rm Supp}({\rm polyhed}(p)({\rm New}(s)))$ and ${\rm New}(c_i)$ are polyhedra for all $i$, this linear functional $\ell$ is minimised on a face of ${\rm polyhed}(p)(({\rm New}(s)))$. For any vertex of this face, $\tilde{\gamma}$ is the unique vertex of ${\rm New}(s)$ associated with it.  Since there are finitely many points of the form $-(v_i-v_j)/(i-j)$, we conclude that ${\rm New}(s)$ has only finitely many vertices.  \end{proof}

The following property starts a link between the classical and the polyhedral settings. Since convex hull and Minkowski sum are well-defined for arbitrary sets, the evaluation $\phi(\mathcal{B})$ of a polyhedral polynomial $\phi:=\oplus_i (Q_i \odot Y^{i})$ at a convex set $\mathcal{B}$ can be defined as $\oplus_i (Q_i \odot \mathcal{B}^{\odot i})$.  Furthermore, since the set of convex sets is closed under both these operations, $\phi(\mathcal{B})$ is a convex set.

\begin{proposition}\label{rootpres_prop}
Let ${\bf s}=(s_1,\dots,s_h) \in \mathbb{K}^h$ be a non-zero root of a polynomial $p \in \mathbb{K}_P[x_1,\dots,x_h]$ and let ${\rm New}({\bf s})=({\rm New}(s_1),\dots,{\rm New}(s_h))$.  Every vertex of ${\rm polyhed}(p)({\rm New}({\bf s}))$ is shared by at least two distinct summands.
\end{proposition}

\begin{proof}
Suppose that $p=\sum_{\alpha|c_{\alpha} \neq 0} c_{\alpha} {\bf x}^{\alpha}$.   Since ${\bf s}$ is a root of $p$, we have $\sum_{\alpha|c_{\alpha} \neq 0} c_{\alpha} {\bf s}^{\alpha}=0$.  
Since ${\bf s}$ is non-zero, $\sum_{\alpha|c_{\alpha} \neq 0} c_{\alpha} {\bf s}^{\alpha}$ has a non-zero summand. Consider a vertex of $\oplus_{\alpha|c_{\alpha} \neq 0} {\rm New}(c_{\alpha} {\bf s}^{\alpha})$. It must be a vertex of some summand ${\rm New}(c_{\tilde{\alpha}} {\bf s}^{\tilde{\alpha}})$,  say where $\tilde{\alpha}=(\tilde{\alpha}_1,\dots,\tilde{\alpha}_h)$ and hence, is of the form  $u+\sum_{i=1}^{h} \tilde{\alpha}_i \cdot v_i$  where $u$ and $v_i$ are vertices of ${\rm New}(c_{\tilde{\alpha}})$ and ${\rm New}(s_i)$, respectively. Since the coefficient of ${\bf t}^{u+\sum_{i=1}^{h} \tilde{\alpha}_i \cdot v_i}$ in (the $\omega$-positive Puiseux series) $c_{\tilde{\alpha}} {\bf s}^{\tilde{\alpha}}$ is non-zero and  $\sum_{\alpha|c_{\alpha} \neq 0} c_{\alpha} {\bf s}^{\alpha}=0$, there must be another summand $c_{\hat{\alpha}} {\bf s}^{\hat{\alpha}}$ that contains $u+\sum_{i=1}^{h} \tilde{\alpha}_i \cdot v_i$ in its support. Furthermore, $u+\sum_{i=1}^{h} \tilde{\alpha}_i \cdot v_i$ is a vertex of ${\rm New}(c_{\hat{\alpha}} {\bf s}^{\hat{\alpha}})$ (since it is a vertex of  $\oplus_{\alpha|c_{\alpha} \neq 0} {\rm New}(c_{\alpha} {\bf s}^{\alpha})$). Hence, every vertex of $\oplus_{\alpha|c_{\alpha} \neq 0} {\rm New}(c_{\alpha} {\bf s}^{\alpha})$ is shared by at least two distinct summands. Since ${\rm New}(c_{\alpha} {\bf s}^{\alpha})={\rm New}(c_{\alpha}) \odot ({\rm New}({\bf s}))^{\odot \alpha}$  for each $\alpha$, the corresponding summands of ${\rm polyhed}(p)({\rm New}({\bf s}))$  and $\oplus_{\alpha|c_{\alpha} \neq 0} {\rm New}(c_{\alpha} {\bf s}^{\alpha})$ are equal. This completes the proof of the proposition. \end{proof}

Note that Proposition \ref{rootpres_prop} is vacuous if ${\rm polyhed}(p)({\rm New}({\bf s}))=0_{\mathcal{P}_{\omega}}$.
In the following proposition, we show that there are roots of ${\rm polyhed}(p)$ that are arbitrarily close to ${\rm New}(s)$. This justifies studying (polyhedral) solutions to ${\rm polyhed}(p)$ in order to obtain information about solutions to $p$.

\begin{proposition}\label{polyhed_prop}
With the same setting as in Proposition \ref{finver_prop}, $V({\rm New}(s)) \odot C$ (where $V(.)$ is the vertex set) is a root of ${\rm polyhed}(p)$ for any polyhedral cone $C$ that contains the recession cone of ${\rm New}(s)$.
\end{proposition}

\begin{proof}
Let $\phi={\rm polyhed}(p)$, $B={\rm New}({\bf s})$ and $P=V({\rm New}({\bf s})) \odot C$. We have $\phi(P)=\phi(B) \odot C$ and the $i$-th summand $\Psi_i(P)=\Psi_i(B) \odot C$ for each $i \in {\rm Supp}(\phi)$.  For a polyhedron $\tilde{P}$ and a cone $K$, we have $V(\tilde{P} \odot K) \subseteq V(\tilde{P})$. Furthermore, a vertex $u$ of $\tilde{P}$ is also a vertex of $\tilde{P} \odot K$ if and only if $K^{\star} \cap N_{\tilde{P}}(u)$ is full-dimensional. Hence, $V(\phi(P)) \subseteq V(\phi(B))$ and $V(\Psi_i(P)) \subseteq V(\Psi_i(B))$ for each $i$. 

Consider $v \in V(\phi(B))$. By construction, $N_{\phi(B)}(v) \subseteq N_{\Psi_i(B)}(v)$ for every $i \in {\rm Supp}(\phi)$ such that $v \in V(\Psi_i(B))$. We deduce that if $v \in V(\phi(P))$, then $v \in V(\Psi_i(P))$ for all $i$ such that $v \in V(\Psi_i(B))$.  By Proposition \ref{rootpres_prop}, every vertex of $\phi(P)$ is also shared by at least two of its summands. 
\end{proof}




\section{The Local Case}


In this section, we formalise the notion of local compatible system that plays a key role in Theorem \ref{minweychar_theointro}.  We start by defining the notion of a generic polyhedral polynomial.


Recall that the polynomial $\phi$ is of the form $Q_{i_{r+1}} \odot Y^{i_{r+1}} \oplus Q_{i_r} \odot  Y^{i_r} \oplus \dots \oplus Q_{i_0} \odot Y^{i_0}$ where $0 \leq i_0<i_1<i_2<\dots<i_{r+1}=d$. We refer to the set $\{i_{r+1},\dots,i_0\}$ as the support of $\phi$ and denote it by ${\rm Supp}(\phi)$. Let $M$ be the Minkowski sum $Q_{i_0} \odot \cdots \odot Q_{i_{r+1}}$ of the coefficients of $\phi$.  Any vertex $\nu$ of $M$ can be uniquely expressed, via its Minkowski decomposition, as $\nu_{i_0}+\dots+\nu_{i_{r+1}}$ where $\nu_{i_k}$ is a vertex of $Q_{i_k}$ for $k$ from $0$ to $r+1$ \cite[Page 16]{Wei07}.  Let $\rho_{i,j}^{(\nu)}:=-(\nu_i-\nu_j)/(i-j)$ for distinct $i,~j \in 
{\rm Supp}(\phi)$.

\begin{definition}{\rm{(\bf Genericity})}\label{generpoly_def}
A polyhedral polynomial $\phi$ is called generic if for every vertex $\nu$ of $M$, the points $\rho_{\hat{i},\hat{j}}^{(\nu)}=\rho_{\tilde{i},\tilde{j}}^{(\nu)}$ implies that $\{\hat{i},\hat{j}\}= \{\tilde{i},\tilde{j}\}$.
\end{definition}

In other words, the genericity condition (Definition \ref{generpoly_def}) asserts that for any fixed vertex $\nu$ of $M$, the points $\rho_{i,j}^{(\nu)}$ are all distinct.   The following proposition justifies the terminology ``generic'' in Definition \ref{generpoly_def}. Given any polyhedral polynomial $\phi$ and an $\epsilon_0>0$, let $A(\epsilon_0)$ be the set of elements $\{\zeta_{\gamma,i}\}$ where each $\zeta_{\gamma,i} \in  \mathbb{S}_{\gamma,i}^n$ (a copy of the unit sphere in $\mathbb{R}^n$), ${\gamma \in V(Q_i),i \in {\rm Supp}(\phi)}$ such that for every $0<\epsilon \leq \epsilon_0$ the following properties hold:

\begin{itemize}

\item For every $i \in {\rm Supp}(\phi)$ and $\gamma \in V(Q_i)$, the point $\gamma^{(\epsilon)}=\gamma+\epsilon \cdot \zeta_{\gamma,i}$ is a vertex of $Q^{(\epsilon)}_i$ where $Q^{(\epsilon)}_i={\rm CH}(\{\gamma^{(\epsilon)}\}_{\gamma \in V(Q_i)}) \odot {\rm MWC}(Q_i)$ (recall that ${\rm CH(.)}$ is the convex hull and ${\rm MWC}(.)$ is the Minkowski-Weyl cone).

\item The polyhedral polynomial $\phi^{(\epsilon)}:= \oplus_{i \in {\rm Supp}(\phi)} (Q^{(\epsilon)}_i \odot Y^i)$ is generic. 
 
 \end{itemize}

\begin{proposition}
There is an $\epsilon_0>0$ such that the set $A(\epsilon_0)$ has full measure  as a subset of the product of the spheres $\mathbb{S}_{\gamma,i}^n$. 
\end{proposition}
\begin{proof} 
By a continuity argument, there is an $\epsilon_0>0$ such that  $\gamma^{(\epsilon)}=\gamma+\epsilon \cdot \zeta_{\gamma,i}$ is a vertex of $Q^{(\epsilon)}_i$ for every $i \in {\rm Supp}(\phi)$ and $\gamma \in V(Q_i)$.  For any quadruple $\gamma_i^{(\epsilon)},\gamma_j^{(\epsilon)},\gamma_{\tilde{i}}^{(\epsilon)},\gamma_{\tilde{j}}^{(\epsilon)}$ such that $i \neq j,~\tilde{i} \neq \tilde{j}$ and $\gamma^{(\epsilon)}_k$ is a vertex of $Q^{(\epsilon)}_k$ for all valid $k$, the condition $(\gamma_i^{(\epsilon)}-\gamma_j^{(\epsilon)})/(i-j)=(\gamma_{\tilde{i}}^{(\epsilon)}-\gamma_{\tilde{j}}^{(\epsilon)})/(\tilde{i}-\tilde{j})$ imposes a non-trivial affine condition on $\zeta_{\gamma_i,i}, \zeta_{\gamma_j,j}, \zeta_{\gamma_{\tilde{i}},\tilde{i}}, \zeta_{\gamma_{\tilde{j}},\tilde{j}}$. The statement then follows from the observation that there are only finitely many such affine conditions.  
\end{proof}

\subsection{Local Compatible Systems}\label{lcs_sect}


We start by constructing certain labelled fans from the polyhedral polynomial.


{\bf Labelled Fans of a Polyhedral Polynomial:}  We start by defining two families of hyperplanes.  Fix a vertex $v$ of $M$ and let $\sum_{i} v_i$ be its Minkowski decomposition, define $H^{(1,v)}_{i,j,k}:= \{\mu \in (\mathbb{R}^n)^{\star}  |~\mu(v_k+k  \cdot \rho_{i,j}^{(v)})=\mu(v_i+i \cdot \rho_{i,j}^{(v)}) \}$ and $H^{(2,v)}_{i,j,\tilde{i},\tilde{j}}:=\{\mu \in (\mathbb{R}^n)^{\star}|~\mu(\rho_{i,j}^{(v)})=\mu(\rho_{\tilde{i},\tilde{j}}^{(v)})\}$ for $i, j,\tilde{i}, \tilde{j},k \in {\rm Supp}(\phi)$ such that $i \neq j$ and $\tilde{i} \neq \tilde{j}$. Observe that $v_i+i \cdot \rho_{i,j}^{(v)} = v_j +j \cdot \rho_{i,j}^{(v)}=(i \cdot v_j-j \cdot v_i)/(i-j)$. Also, note that $H^{(1,v)}_{i,j,k}$ and $H^{(2,v)}_{i,j,\tilde{i},\tilde{j}}$ are hyperplanes when $v_k+k  \cdot \rho_{i,j}^{(v)} \neq v_i+i \cdot \rho_{i,j}^{(v)}$ and $\rho_{i,j}^{(v)} \neq \rho_{\tilde{i},\tilde{j}}^{(v)}$, respectively. In this case, we refer to them as \emph{non-degenerate}. Let $\mathcal{ND}_v$ denote the set of non-degenerate hyperplanes of the form either $H^{(1,v)}_{i,j,k}$ or $H^{(2,v)}_{i,j,\tilde{i},\tilde{j}}$.  In the case where $\mathcal{ND}_v$ is non-empty,  let $\mathcal{H}_v$ be the hyperplane arrangement $\cup_{H \in \mathcal{ND}_v} H$ in $(\mathbb{R}^n)^{\star}$. Note that $\mathcal{H}_v$ is a central hyperplane arrangement.  It is the (inner or outer) normal fan of the zonotope $Z_v$ given by the Minkowski sum of a choice of normal vectors one for each hyperplane \cite[Chapter 7.3]{Zie07}.  
Consider the restriction $\mathcal{H}_v| N_M(v)$ of $\mathcal{H}_v$ to the inner normal cone $N_M(v)$ of $v$ (as a vertex of $M$). In other words, $\mathcal{H}_v| N_M(v)$ is the inner normal fan of the Minkowski sum of $Z_v \odot (N_M(v))^{\star}$ of $Z_v$ with the tangent cone of $v$. It is a polyhedral fan of dimension $n$ in $(\mathbb{R}^n)^{\star}$.  We denote this fan by $\mathcal{F}_v$. 


We label each maximal cone of  $\mathcal{F}_v$ with some additional data. The labels are of two types: primary and secondary. A maximal cone $C$ of  $\mathcal{F}_v$ has $(i,j)$ as a \emph{primary label} if $\mu(v_i+ i \cdot \rho_{i,j}^{(v)}) (=\mu(v_j+j \cdot \rho_{i,j}^{(v)}))  \leq \mu(v_k+k \cdot \rho_{i,j}^{(v)})$  for all $k \in \{i_0,\dots,i_{r+1}\}$ and all $\mu \in C$. The cone $C$ has $(i,j,\tilde{i},\tilde{j})$ as a \emph{secondary label} if  $\mu(\rho_{i,j}^{(v)}) \leq \mu(\rho_{\tilde{i},\tilde{j}}^{( v)})$ for all $\mu \in C$.  Equivalently,  $(i,j,\tilde{i},\tilde{j})$ is a secondary label of $C$ if $\rho_{\tilde{i},\tilde{j}}^{(v)}-\rho_{i,j}^{(v)}$ is contained in the dual cone $C^{\star}$ of $C$. We refer to the corresponding inequalities as primary and secondary inequalities, respectively.  Note that, in general, a maximal cone can have multiple primary as well as secondary labels. We refer to Section \ref{example_sect} for examples.

\begin{remark}\label{primlab_rem}
\rm
An equivalent interpretation of primary labels is the following.  The maximal cone $C \in \mathcal{F}_v$ has $(i,j)$ as a primary label if for all $\mu \in C$ and all $k \in {\rm Supp}(\phi) \setminus \{i\}$ the following inequalities hold:
$\mu(\rho_{i,j}^{(v)}) \geq \mu(\rho_{i,k}^{(v)})$ if $k>i$ and $\mu(\rho_{i,j}^{(v)}) \leq  \mu(\rho_{i,k}^{(v)})$, otherwise. \qed
\end{remark}


\begin{proposition} {\rm{(\bf Algebraically Closed Property of $\mathcal{P}_\omega$)}}
Every polyhedral polynomial $\phi \in \mathcal{A}[Y]$ has an affine cone solution. 
\end{proposition}
\begin{proof}
Let $v$  be the unique vertex that minimises the linear functional $\omega$ over $M$. Consider the convex hull  of the set $\{\rho^{(v)}_{i,j}\}_{i,j}$ and suppose that $\rho^{(v)}_{\hat{i},\hat{j}}$ is the unique vertex that minimises $\omega$ over it.  The affine cone $\rho^{(v)}_{\hat{i},\hat{j}} \odot (N_M(v))^{\star}$ is a solution to $\phi$ and follows from Remark \ref{primlab_rem}.
\end{proof}



\begin{definition} {\rm ({\bf Local Compatible Systems})} \label{lcs_def}
Fix a vertex $v$ of $M$.  A local compatible system of $\phi$ (at $v$) is a pair $(\mathcal{C},\mathcal{I})$ where $\mathcal{C}=\{C_1,\dots,C_t\}$ is a collection of maximal cones  of $\mathcal{F}_v$ and $\mathcal{I}=\{(i_1,j_1),\dots,(i_t,j_t)\}$ is a corresponding collection of pairs of distinct integers $i_k,j_k \in {\rm Supp}(\phi)$ such that the following conditions are satisfied:

\begin{enumerate}

\item \label{comp0_con} Let $\mathfrak{T}:=\{ \rho^{(v)}_{i_k,j_k}\}_{k=1}^{t}$. For each $\gamma \in \mathfrak{T}$,  the set $C_{\gamma}:=\cup_{\{k|~ \rho^{(v)}_{i_k,j_k}=\gamma\}} C_k$  is a convex polyhedral cone. 

\item \label{comp1_con} For each $k$, the pair $(i_k,j_k)$ is a primary label of $C_k$.

\item \label{comp2_con} For every pair of distinct $\alpha,\beta \in \{1,\dots,t\}$, the tuple $(i_{\alpha},j_{\alpha},i_{\beta},j_{\beta})$ is a secondary label of $C_{\alpha}$.

Suppose that $\mathfrak{f}$ is a facet of a unique cone $C_l$ among $C_1,\dots,C_t$ and suppose that there is a maximal cone $C$, say that is not in $\mathcal{C}$ that shares this facet \footnote{Note that such a maximal cone is unique, if it exists. 
}.

\item \label{comp4_con}   Either $C$  must  not have the primary label  $(i_l,j_l)$, otherwise there must an index $\beta \in \{1,\dots,t\}$ such that $(i_l,j_l,i_{\beta},j_{\beta})$ is not a secondary label of $C$. 

\end{enumerate}

\end{definition}

\begin{remark}\label{supplcs_rem}\rm
The collection $\mathcal{C}$ is ``almost'' the normal fan of a polyhedron except that the union of its elements is not necessarily convex. This union $\cup_{k=1}^{t}C_k$ is called the \emph{support} of the local compatible system. \qed
\end{remark}


\subsection{Vertex-Cone Collections}

In this section, we develop the notion of vertex-cone collection with a view towards Theorem \ref{genminwey_theo}.  Let $W$ be a finite dimensional real vector space.  

\begin{definition}{\rm ({\bf  Vertex-Cone Collection})}\label{vercon-def}
A collection $\mathfrak{G}=\{(v,N(v))\}_{v \in G}$ where $G$ is a non-empty, finite subset  of $W$ and each $N(v)$ is a full-dimensional convex polyhedral cone in $W^{\star}$ is called a vertex-cone collection (VCC) if  for every pair of distinct points $v,u \in G$ and for every $\ell \in N(v)$, the inequality $\ell(v) \leq \ell(u)$ holds.\end{definition}


Note that the inequality in Definition \ref{vercon-def} ensures that $N(u)$ and $N(v)$ do not intersect in their interior.


The elements of $G$ are called the \emph{vertices} of $\mathfrak{G}$ and the cone $N(v)$ for $v \in G$ in Definition \ref{vercon-def} is called the \emph{normal cone} of $v$. We also denote it by $N_G(v)$ or $N_{\mathfrak{G}}(v)$ as appropriate.
 We refer to $\cup_{v \in G} N(v)$ as the \emph{support} of $\mathfrak{G}$ and denote it by ${\rm Supp}(\mathfrak{G})$. The support of $\mathfrak{G}$ is the analogue of the normal cone of a polyhedron.   
 
 \begin{remark} \label{polyvcc_rem} \rm
 As mentioned in the introduction, a vertex-cone collection only captures the behaviour of the vertices of a polyhedron with respect their normal cones. More precisely, for any polyhedron $P$ in $W$, the collection $\{(v,N(v))\}_{v \in V(P)}$  where $V(P)$ is the vertex set of $P$ and $N(v)$ is the normal cone of $P$ is a vertex-cone collection. This is called the \emph{associated vertex-cone collection}  of $P$. The support of such a vertex-cone collection is a convex cone whereas the support of an arbitrary vertex-cone collection need not be convex.  \qed
\end{remark}

\subsubsection{Operations on Vertex-Cone Collections}

Our primary goal here is to introduce the notion of a VCC solution to a polyhedral polynomial. En route to this goal, we define basic operations on vertex-cone collections  analogous to convex hull and Minkowski sum. 
Let $\mathfrak{G_1}=\{(v,N_{G_1}(v))\}_{v \in G_1}$ and $\mathfrak{G}_2=\{(v,N_{G_2}(v))\}_{v \in G_2}$ be vertex-cone collections. 

\begin{definition}{\rm{({\bf Convex Hull of Vertex-Cone Collections})}}\label{convhullvcc_def}
Suppose that the supports of $\mathfrak{G}_1$ and $\mathfrak{G}_2$ are equal. The convex hull  $\mathfrak{G_1} \oplus \mathfrak{G_2}$ of  $\mathfrak{G}_1$ and  $\mathfrak{G}_2$  is the collection $\{(v,N_{\hat{G}}(v))\}_{v \in \hat{G}}$ where $\hat{G}$ is defined as follows.  It consists of all elements in $G_1 \cap G_2$ such that $N_{G_1}(v) \cap N_{G_2}(v)$ is a full-dimensional cone along with elements $v$ where $v \in G_i \setminus G_j$ for distinct $i,j \in \{1,2\}$ such that there is an element $\ell$ in the interior of $N_{G_i}(v)$ such that $\ell(v)<\ell(u)$ for all $u \in G_j$. The corresponding normal cones are as follows: $N_{\hat{G}}(v):=N_{G_1}(v) \cap N_{G_2}(v)$ if $v \in G_1 \cap G_2$ and $N_{\hat{G}}(v):=\{ \ell \in N_{G_i}(v)|~\ell(v) \leq \ell(u)~\forall u \in G_j\}$ if $v \in G_i \setminus G_j$. 
\end{definition} 

Note that the normal cones $N_{\hat{G}}(v)$ in Definition \ref{convhullvcc_def} are full-dimensional convex polyhedral cones. We emphasise that the convex hull of vertex-cone collections (Definition \ref{convhullvcc_def}) requires the supports to be equal. Next, we turn to Minkowski sums of vertex-cone collections. We start with a version of Minkowski decomposition for vertex-cone collections. 

\begin{proposition}{\rm{({\bf Minkowski Decomposition})}}\label{mindecgenpol_prop}
Fix a positive integer $k$ and let $\mathfrak{G}_i=\{(v,N_{G_i}(v))\}_{v \in G_i}$ be a vertex-cone collection for each $i$ from one to $k$.  If $\cap_{i=1}^{k} N_{G_i}(v_i)$ is a full-dimensional cone for vertices $v_i \in G_i$, then the equation ${\bf x}_1+\dots+{\bf x}_k=v_1+\dots+v_k$ in variables ${\bf x}_i \in G_i$  has a unique solution, namely ${\bf x_i}=v_i$. 
\end{proposition}
The proof of Proposition \ref{mindecgenpol_prop} is relatively straightforward and is hence, omitted. As before, we refer to $v_1+\dots+v_k$ as the {\emph{Minkowski decomposition}} of this vector. 


\begin{definition}{\rm{({\bf Minkowski Sum of Vertex-Cone Collections})}}
Suppose that the intersection of the supports of $\mathfrak{G}_1$ and $\mathfrak{G}_2$ is full-dimensional.  The Minkowski sum $\mathfrak{G_1} \odot \mathfrak{G_2}$ of $\mathfrak{G_1}$ and $\mathfrak{G_2}$ is defined as the collection $\{(v,N_{\tilde{G}}(v))\}_{v \in \tilde{G}}$ where $\tilde{G}=\{v_1+v_2|~v_1 \in G_1,v_2 \in G_2 \text{ and } N_{G_1}(v_1) \cap N_{G_2}(v_2) \text{~is a full-dimensional cone}\}$ and for each $w \in \tilde{G}$, the cone $N_{\tilde{G}}(w)$ is defined as $N_{G_1}(v_1) \cap N_{G_2}(v_2)$ where $v_1+v_2$ is the Minkowski decomposition of $w$.
\end{definition}

 Both these definitions were derived by examining normal fans of polyhedra under the corresponding operations.   Observe that the definition of convex hull of vertex-cone collections is somewhat more complicated than that of Minkowski sum. This is rooted in the fact that the behaviour  of normal fans of polyhedra under convex hull is more complicated than under Minkowski sum.   The following proposition shows that the result of the two operations are indeed vertex-cone collections. Suppose that $\mathfrak{G}_1=\{(v,N_{G_1}(v))\}_{v \in G_1}$ and $\mathfrak{G}_2=\{(v,N_{G_2}(v))\}_{v \in G_2}$ are vertex-cone collections. 
 
 \begin{proposition}\label{vccop_prop}  If  ${\rm Supp}(\mathfrak{G}_1)={\rm Supp}(\mathfrak{G}_2)$, then the convex hull $\mathfrak{G}_1 \oplus \mathfrak{G}_2$ is a vertex-cone collection. 
If   ${\rm Supp}(\mathfrak{G}_1) \cap {\rm Supp}(\mathfrak{G}_2)$ is full-dimensional, then the Minkowski sum $\mathfrak{G}_1 \odot \mathfrak{G}_2$ is a vertex-cone collection. 
 \end{proposition}

 The following lemma relates the supports of $\mathfrak{G}_1 \oplus \mathfrak{G}_2$ and $\mathfrak{G}_1 \odot \mathfrak{G}_2$ to its constituents $\mathfrak{G}_1$ and $\mathfrak{G}_2$.
 
 \begin{lemma}\label{supp_lemma}
 With the same hypothesis on the supports as in Proposition \ref{vccop_prop},  ${\rm Supp}(\mathfrak{G}_1 \oplus \mathfrak{G}_2)={\rm Supp}(\mathfrak{G}_1)={\rm Supp}(\mathfrak{G}_2)$  and ${\rm Supp}(\mathfrak{G}_1 \odot \mathfrak{G}_2)={\rm Supp}(\mathfrak{G_1}) \cap {\rm Supp}(\mathfrak{G_2})$.
 \end{lemma}



 As in the case of polyhedra, both the operations $\oplus$ and $\odot$ are commutative and associative.  We say that collections $\mathfrak{G}_1=\{(v,N_{G_1}(v))\}_{v \in G_1}$ and $\mathfrak{G}_2=\{(v,N_{G_2}(v))\}_{v \in G_2}$ are \emph{equal} if there is a bijection  $\psi$ between $G_1$ and $G_2$ such that  $\psi(u)=u$ and $N_{G_1}(u)=N_{G_2}(\psi(u))$ for every  $u \in G_1$.
 
\begin{proposition}\label{commassogp_prop}
 Let $\mathfrak{G}_i=\{(v,N_{G_i}(v))\}_{v \in G_i}$ be vertex-cone collections. With appropriate hypothesis on the supports, the following equalities hold:
 \begin{itemize}
 \item ({\bf Commutativity}) The vertex-cone collections $\mathfrak{G}_i \oplus \mathfrak{G}_j=\mathfrak{G}_j \oplus \mathfrak{G}_i$ and $\mathfrak{G}_i \odot \mathfrak{G}_j=\mathfrak{G}_j \odot \mathfrak{G}_i$.
 \item  ({\bf Associativity for Convex Hull})  The  collections $\mathfrak{G}_1 \oplus (\mathfrak{G}_2 \oplus \mathfrak{G}_3)$ and $(\mathfrak{G}_1 \oplus \mathfrak{G}_2) \oplus \mathfrak{G}_3$ are both well-defined vertex-cone collections and are equal. 
 
\item   ({\bf Associativity for Minkowski Sum}) The collections  $\mathfrak{G}_1 \odot (\mathfrak{G}_2 \odot \mathfrak{G}_3)$ and $(\mathfrak{G}_1 \odot \mathfrak{G}_2) \odot \mathfrak{G}_3$ are both  well-defined vertex-cone collections and are equal.
  \end{itemize}
 \end{proposition}
 
 Distributivity comes next in this sequence that we do not treat here.  Since the proofs of Propositions \ref{vccop_prop} and \ref{commassogp_prop}, and Lemma \ref{supp_lemma} are not central to  the rest of the article, they have been included in an appendix (Appendix \ref{vcc_app}). 
 
 
 \subsubsection{Polynomials and their Roots}

 
The associativity part of Proposition \ref{commassogp_prop} implies that finite sums and finite products, with appropriate hypotheses on the supports, of vertex-cone collections are well-defined. This allows us to evaluate polyhedral polynomials at vertex-cone collections. Fix a positive integer $d$, integers $1 \leq  r \leq d-1$ and $1 \leq  i_0 < \dots < i_{r}<i_{r+1}=d$.  Let $\phi=Q_d \odot Y^d \oplus Q_{i_r} \odot Y^{i_r} \oplus \dots \oplus Q_{i_0} \odot Y^{i_0}$ be a polyhedral polynomial. Note that, since $i_0 \geq 1$, the polynomial $\phi$ does not have a constant term. In the following, we denote the inner normal cone of a polyhedron $P$ by $N(P)$. For a coefficient polyhedron $Q_k$ of $\phi$, we denote its associated vertex-cone collection (recall Remark \ref{polyvcc_rem}) by $\mathfrak{Q}_k$.



\begin{prop-def}{\rm ({\bf Evaluation at a VCC})}\label{genpolyneva_def}
   Given a vertex-cone collection  $\mathfrak{P}$ such that ${\rm Supp}(\mathfrak{P}) \subseteq \cap_{j=1}^{r+1}N(Q_j)$, the expression $\mathfrak{Q}_{i_1} \odot \mathfrak{P}^{\odot i_1} \oplus \dots \oplus \mathfrak{Q}_{i_{r+1}} \odot \mathfrak{P}^{\odot i_{r+1}}$ is  well-defined and is a vertex-cone collection with support ${\rm Supp}(\mathfrak{P})$,  referred to as the evaluation $\phi(\mathfrak{P})$ of $\phi$ at $\mathfrak{P}$. 
\end{prop-def}

The well-definedness in Proposition-Definition \ref{genpolyneva_def} can be seen as follows. By Lemma \ref{supp_lemma}, the support hypotheses required for convex hull and Minkowski sum of VCC are satisfied. By Proposition \ref{commassogp_prop}, the expression does not require specifying the order in which the multiplications in each term and the addition of the terms are to be performed.  The claim that $\phi(\mathfrak{P})$ is a VCC with support ${\rm Supp}(\mathfrak{P})$ follows from Proposition \ref{vccop_prop} and Lemma \ref{supp_lemma}.


 The term  $\mathfrak{Q}_i \odot \mathfrak{P}^{\odot i}$ is called the $i$-th summand of $\phi(\mathfrak{P})$. Next, we define the notion of a root of a polynomial in the setting of vertex-cone collections.



\begin{definition}{\rm (\bf Root of a Polynomial)}\label{vccroot_def}
 A vertex-cone collection $\mathfrak{P}_0$ such that ${\rm Supp}(\mathfrak{P}_0) \subseteq \cap_{k=1}^{r+1} N(Q_k)$ is said to be a root of $\phi$ if every vertex of $\phi(\mathfrak{P}_0)$ is shared by at least two summands, i.e. for every vertex $v$ of  $\phi(\mathfrak{P}_0)$, there exist distinct integers $i,j \in {\rm Supp}(\phi)$ such that $v$ is a vertex of both $\mathfrak{Q}_i \odot \mathfrak{P}_0^{\odot i}$ and $\mathfrak{Q}_j \odot \mathfrak{P}_0^{\odot j}$. 
\end{definition}

\begin{remark}
\rm  The following remarks are called for:

\begin{itemize}

\item The conditions on the support of $\mathfrak{P}$ and $\mathfrak{P}_0$ (in \ref{genpolyneva_def} and \ref{vccroot_def}, respectively)  ensure that the resulting convex hull operation is well-defined. They can be eliminated but this leads to more complicated definitions, so we stick to this version.

\item These definitions (\ref{genpolyneva_def} and  \ref{vccroot_def}) require the polynomial $\phi$ to not have the constant term, i.e. $0 \notin {\rm Supp}(\phi)$. This restriction has been imposed to ensure that each summand of $\phi(\mathfrak{P}_0)$ has the same support: a condition needed for the convex hull operation (Definition \ref{convhullvcc_def}). However, for any polynomial $\tilde{\phi}=\oplus_{i} (Q_i \odot Y^{\odot i} )  \in \mathcal{A}[Y]$, the polynomial $\hat{\phi}=\oplus_{i} (Q_i \odot Y^{\odot (i+1)})$ is of the required form.    The set of roots of $\hat{\phi}$ is precisely the set of roots of $\tilde{\phi}$ along with the element $0_{\mathcal{A}}$. 
\end{itemize} \qed
\end{remark}

Analogous to the case of polyhedra, we identify certain extremal solutions: the Minkowski-Weyl minimal ones. 


\begin{definition}{\rm (\bf Minkowski-Weyl Minimality)}\label{minweylgenpoly_def}
Let $V$ be a finite subset of $W$.  A vertex-cone collection $\mathfrak{P}_0$ that is a root of  $\phi$ is said to be a Minkowski-Weyl minimal solution to $\phi$ with respect to $V$ if it has vertex set $V$ and has maximal support among all  vertex-cone collections with the same vertex set that are also roots of $\phi$.
\end{definition}

In this case, we say that $\mathfrak{P}_0$ is a Minkowski-Weyl minimal root to mean that its support is maximal among all vertex-cone collection roots of  $\phi$ with the same vertex set as $\mathfrak{P}_0$.

Note that Definition \ref{minweylgenpoly_def} is a dual version of the corresponding definition for polyhedra. Recall that if  cones $C_1 \subseteq C_2$, then $C_2^\star \subseteq C_1^\star$. The reversal of minimality to maximality arises from the fact that the Minkowski-Weyl cone of a polyhedron is dual to its normal cone and that the support of a vertex-cone collection plays the role of its normal cone. 


In the following, we characterise Minkowski-Weyl minimal vertex-cone collection solutions to $\phi$ in the local case, i.e. solutions whose support is contained in the normal fan of a fixed vertex $v$ of $M$ and under the assumption that $\phi$ is a generic polyhedral polynomial.  

Let $\phi_{\rm gen}$ be a generic polyhedral polynomial with $0 \notin {\rm Supp}(\phi_{\rm gen})$. 
Given a local compatible system associated with $\phi_{\rm gen}$, we associate a vertex-cone collection to it and  show (in this generic, local setting) that this induces a bijective correspondence between local compatible systems and Minkowski-Weyl minimal vertex-cone collection solutions. In the following, we associate a vertex-cone collection to a local compatible system of $\phi_{\rm gen}$.  Recall that $\mathfrak{T}=\{ \rho^{(v)}_{i_k,j_k}\}_{k=1}^{t}$ and the set $C_{\gamma}:=\cup_{\{k|~ \rho^{(v)}_{i_k,j_k}=\gamma\}} C_k$.

\begin{proposition}
Let $(\mathcal{C},\mathcal{I})$ be a local compatible system of $\phi_{\rm gen}$. The set $\{(\gamma,C_{\gamma})\}_{\gamma \in \mathfrak{T}}$ is a vertex-cone collection.
\end{proposition} 
\begin{proof}
By Condition \ref{comp0_con}, the set $C_{\gamma}:=\cup_{\{k|~ \rho^{(v)}_{i_k,j_k}=\gamma\}} C_k$ is a convex polyhedral cone for each $\gamma \in \mathfrak{T}$.  By Condition \ref{comp2_con} (in the definition of local compatible system), the inequality $\sigma(\rho^{(v)}_{i_k,j_k}) \leq \sigma(\rho^{(v)}_{i_l,j_l})$ holds for every $k,l \in [1,\dots,t]$ and for all $\sigma \in C_k$. This implies that $\sigma(\gamma) \leq \sigma(\tilde{\gamma})$ for all $\gamma,\tilde{\gamma} \in \mathfrak{T}$ and for all $\sigma \in C_{\gamma}$.  Hence, $\{(\gamma,C_{\gamma})\}_{\gamma \in \mathfrak{T}}$ is a vertex-cone collection.  
\end{proof}

We refer to this collection $\{(\gamma,C_{\gamma})\}_{\gamma \in \mathfrak{T}}$ as the \emph{associated vertex-cone collection} of $(\mathcal{C},\mathcal{I})$. 

\begin{theorem}\label{genminwey_theo} {\rm ({\bf Generic-Local VCC Case})}
Let $\phi_{\rm gen}$  be a generic polyhedral polynomial with $0 \notin {\rm Supp}(\phi_{\rm gen})$. Fix a vertex $v$ of $M$.  A vertex-cone collection whose support is contained in $N_M(v)$ is a Minkowski-Weyl minimal solution to $\phi_{\rm gen}$ if and only if it is the associated vertex-cone collection of a local compatible system $(\mathcal{C},\mathcal{I})$ of $\phi_{\rm gen}$ such that $\cup_{k=1}^{t}C_k \subseteq N_M(v)$.\end{theorem}


Before the proof of Theorem \ref{genminwey_theo}, we state two core principles that will be a basis for subsequent arguments.

\begin{principle} {\rm ({\bf First Continuation Principle})}\label{firstcont_princ}
If cones $C_1,\dots,C_k$ in $W^{\star}$ satisfy a collection of linear inequalities, then their Minkowski sum $C_1 \odot \dots \odot C_k$ also satisfies this collection of inequalities. 

\end{principle}

\begin{principle}{\rm ({\bf Second Continuation Principle})}\label{seccont_princ}
If a collection of primary and secondary inequalities are valid on a full-dimensional subcone of a maximal cone of $\mathcal{F}_v$ (for a fixed vertex $v$ of $M$), then they are valid on the entire maximal cone. 
\end{principle}

The first continuation principle follows from the definition of Minkowski sum and the linearity of the inequalities. The second continuation principle is a consequence of the definition of $\mathcal{F}_v$.  

\subsection{Proof of Theorem \ref{genminwey_theo}}

Suppose that a vertex-cone collection $\mathfrak{P}_{\rm asso}$ is the associated vertex-cone collection of a local compatible system of $\phi_{\rm gen}$.  We show that  $\mathfrak{P}_{\rm asso}$ is a Minkowski-Weyl minimal solution to $\phi_{\rm gen}$. The following proposition characterising the vertex set and the corresponding normal cones of $\phi_{\rm gen}(\mathfrak{P}_{\rm asso})$ is handy.  We denote the point $i_k \cdot \rho^{(v)}_{i_k,j_k}+v_{i_k}=j_k \cdot  \rho^{(v)}_{i_k,j_k}+v_{j_k}=(i_k \cdot v_{j_k}-j_k \cdot v_{i_k})/(i_k-j_k)$ by $w_{k,i_k,j_k}$. 

 \begin{proposition}\label{phipasso_prop}  The vertices of  $\phi_{\rm gen}(\mathfrak{P}_{\rm asso})$ are precisely of the form  $w_{k,i_k,j_k}$ where $k$ varies from $1$ to $t$. 
  \end{proposition}

\begin{proof}
By the definition of convex hull and Minkowski sum of vertex-cone collections and since $C_k \subseteq \mathcal{F}_{v}$, we know that the vertices of  $\phi_{\rm gen}(\mathfrak{P}_{\rm asso})$ are of the form $i_l \cdot \rho^{(v)}_{i_k,j_k}+v_{i_l}$.  By Condition \ref{comp1_con} (in the definition of local compatible system),  we have the inequality $\sigma(w_{k,i_k,j_k}) \leq \sigma(i_l \cdot \rho^{(v)}_{i_k,j_k}+v_{i_l})$ for all $\sigma \in C_k$ and for all $l$ from $1$ to $r+1$ (recall that ${\rm Supp}(\phi_{\rm gen})=\{i_1,\dots,i_{r+1}\}$).  Hence, we deduce that  $w_{k,i_k,j_k}$ is a vertex for each $k$ from one to $t$. Furthermore, we conclude that there are no other vertices since the support of  $\phi_{\rm gen}(\mathfrak{P}_{\rm asso})$ is contained in the support of $\mathfrak{P}_{\rm asso}$ and this is $\cup_{l=1}^{t}C_l$.
\end{proof}

\begin{remark}
\rm  The vertices $w_{k,i_k,j_k}$ can be shown to be distinct via Conditions \ref{comp1_con} and \ref{comp2_con}. 
\qed
\end{remark}




\begin{proposition}
The vertex-cone collection $\mathfrak{P}_{\rm asso}$ is a Minkowski-Weyl minimal solution to $\phi_{\rm gen}$.
\end{proposition}

\begin{proof}
By Proposition \ref{phipasso_prop}, the vertices of $\phi_{\rm gen}(\mathfrak{P}_{\rm asso})$ are precisely of the form $w_{k,i_k,j_k}$. Furthermore, since $(i_k,j_k)$ is a primary label of $C_k$, the point $w_{k,i_k,j_k}$ is a vertex of both the $i_k$-th and the $j_k$-th summand. Hence, $\mathfrak{P}_{\rm asso}$ is a solution to $\phi_{\rm gen}$.

Next, we turn to the Minkowski-Weyl minimality of $\mathfrak{P}_{\rm asso}$. Suppose for the sake of contraction that $\mathfrak{P}_{\rm asso}$ is not Minkowski-Weyl minimal. Hence, there is a solution $\tilde{\mathfrak{P}}_{\rm asso}$ whose vertex set is the same (as that of $\mathfrak{P}_{\rm asso}$) and whose support strictly contains the support of $\mathfrak{P}_{\rm asso}$.  Note that both $\mathfrak{P}_{\rm asso}$ and $\tilde{\mathfrak{P}}_{\rm asso}$ are vertex-cone collections. Hence, for any vertex  
$w_{k,i_k,j_k}$ its normal cone with respect to  $\mathfrak{P}_{\rm asso}$ is contained in its normal cone with respect to $\tilde{\mathfrak{P}}_{\rm asso}$.  Thus, there is a vertex $w_{k_b,i_{k_b},j_{k_b}}$ for which this containment is strict. 

Since the local compatible system from which $\mathfrak{P}_{\rm asso}$ is constructed satisfies Condition \ref{comp4_con},  there is a maximal cone $C$ of some vertex $u \in M$ such that 
$C \cap N_{\mathfrak{P}_{\rm asso}}(w_{k_b,i_{k_b},j_{k_b}})$ is not full-dimensional and $N_{\tilde{\mathfrak{P}}_{\rm asso}}(w_{k_b,i_{k_b},j_{k_b}}) \setminus N_{\mathfrak{P}_{\rm asso}}(w_{k_b,i_{k_b},j_{k_b}})$ intersects the interior of $C$. Hence, there is an index $i \in {\rm Supp}(\phi_{\rm gen})$ 
such that $u_i+i \cdot w_{k_b,i_{k_b},j_{k_b}}$ is a vertex of $\phi_{\rm gen}(\tilde{\mathfrak{P}}_{\rm asso})$. But this vertex cannot be shared by another summand of $\phi_{\rm gen}(\tilde{\mathfrak{P}}_{\rm asso})$ since otherwise this implies that $C$ has a primary label of the form $(i,j)$ such that $\rho^{(u)}_{i,j}=w_{k_b,i_{k_b},j_{k_b}}$ and hence, contradicts Condition \ref{comp4_con}.   Hence,  $\tilde{\mathfrak{P}}_{\rm asso}$ is not a solution to $\phi_{\rm gen}$. We conclude that $\mathfrak{P}_{\rm asso}$ is a Minkowski-Weyl minimal solution.
\end{proof}

This completes the proof of the right to left implication of Theorem \ref{genminwey_theo}. For the other direction, suppose that $\mathfrak{P}_0$ is a Minkowski-Weyl minimal VCC solution to $\phi_{\rm gen}$ such that ${\rm Supp}(\mathfrak{P}_0) \subseteq N_M(v)$.

\begin{proposition}\label{minweylstruct_prop}
The vertices of $\mathfrak{P}_0$ are of the form $\rho^{(v)}_{i,j}$ for $i,~j \in {\rm Supp}(\phi_{\rm gen})$. The normal cone of each vertex of $\mathfrak{P}_0$ is a finite union of maximal cones of $\mathcal{F}_v$. 
\end{proposition}

\begin{proof}

By Lemma \ref{supp_lemma}, the support of $\phi_{\rm gen}(\mathfrak{P}_0)$ is equal to the support of $\mathfrak{P_0}$. 

Let $\gamma$ be a vertex of $\mathfrak{P}_0$ and let $\sigma$ be a generic element in its normal cone.  Since $\sigma$ is generic, it is minimised at a unique vertex $w$, say of $\phi_{\rm gen}(\mathfrak{P}_0)$. Similarly, over $M$ it is uniquely minimised at $v$.  Since $\mathfrak{P_0}$ is a solution to $\phi_{\rm gen}$, the vertex $w$ is shared by at least two summands of $\phi_{\rm gen}(\mathfrak{P}_0)$.  Suppose that the $i$-th and the $j$-th summands of $\phi_{\rm gen}(\mathfrak{P}_0)$ (for $i \neq j$) satisfy this property, we have $w=v_i+i \cdot \gamma=v_j+j \cdot \gamma$.  Hence, $\gamma=-(v_i-v_j)/(i-j)=\rho^{(v)}_{i,j}$ as claimed. 

We now turn to the second part of the proposition. Suppose for the sake of contradiction that the normal cone of some vertex $\rho^{(v)}_{i,j}$ of $\mathfrak{P}_0$ is not a finite union of maximal cones of $\mathcal{F}_v$. This implies that the normal cone of the vertex $\rho^{(v)}_{i,j}$ intersects a maximal cone $C$ of $\mathcal{F}_v$ in a full-dimensional cone that is strictly contained in $C$. Consider the collection $\tilde{\mathfrak{P}_0}$ defined as $\mathfrak{P}_0 \setminus \{ (\rho^{(v)}_{i,j},N_{\mathfrak{P_0}}(\rho^{(v)}_{i,j})\} \cup \{(\rho^{(v)}_{i,j}, N_{\mathfrak{P_0}}(\rho^{(v)}_{i,j}) \odot C)\}$.  We start by noting that  $\tilde{\mathfrak{P}_0}$ is a vertex-cone collection.  Since $\mathfrak{P}_0$ is a vertex-cone collection, we know that $\ell(\rho^{(v)}_{i,j}) \leq \ell(\rho^{(v)}_{\tilde{i},\tilde{j}})$ for all $\ell \in N_{\mathfrak{P_0}}(\rho^{(v)}_{i,j})$ and for all vertices $\rho^{(\tilde{v})}_{\tilde{i},\tilde{j}}$ of $\mathfrak{P}_0$ (and hence, of $\tilde{\mathfrak{P}_0}$).  Since $N_{\mathfrak{P_0}}(\rho^{(v)}_{i,j}) \cap C$ is full-dimensional, the second continuation principle (Principle \ref{seccont_princ}) asserts that  $\ell(\rho^{(v)}_{i,j}) \leq \ell(\rho^{(v)}_{\tilde{i},\tilde{j}})$ for all $\ell \in C$. By the first continuation principle (Principle \ref{firstcont_princ}),  $\ell(\rho^{(v)}_{i,j}) \leq \ell(\rho^{(v)}_{\tilde{i},\tilde{j}})$ for all $\ell \in   N_{\mathfrak{P_0}}(\rho^{(v)}_{i,j}) \odot C$. Corresponding inequalities for the other vertices of $\tilde{\mathfrak{P}_0}$  are inherited from $\mathfrak{P}_0$. Hence, we conclude that  $\tilde{\mathfrak{P}_0}$ is a vertex-cone collection. 

\b
Next, we show that $\tilde{\mathfrak{P}_0}$ is also a solution to $\phi_{\rm gen}$. We claim that the set of vertices of $\phi_{\rm gen}(\tilde{\mathfrak{P}_0})$ is the same as that of $\phi_{\rm gen}(\mathfrak{P}_0)$.  By the definition of the convex hull operation, the set of vertices of  $\phi_{\rm gen}(\tilde{\mathfrak{P}_0})$ contains the set of vertices of  $\phi_{\rm gen}(\mathfrak{P}_0)$.  Suppose for the sake of contradiction that $\phi_{\rm gen}(\tilde{\mathfrak{P}_0})$ contains an additional vertex. Since the vertex sets of $\mathfrak{P}_0$ and $\tilde{\mathfrak{P}_0}$ are equal and the corresponding normal fans are also equal except for $\rho^{(v)}_{i,j}$,  this vertex must be of the form $k_0 \cdot \rho^{(v)}_{i,j}+v_{k_0}$ for some $k_0 \notin \{i,j\}$.  Furthermore, there is a point $\ell_0 \in C \odot N_{\mathfrak{P_0}}(\rho^{(v)}_{i,j})$ such that  $\ell_0(k_0 \cdot \rho^{(v)}_{i,j}+v_{k_0})<\ell_0(i \cdot \rho^{(v)}_{i,j}+v_i)$. On the other hand, 
$N_{\mathfrak{P_0}}(\rho^{(v)}_{i,j}) \subseteq \mathcal{F}_v$ and since $\phi_{\rm gen}$ is generic, $\rho^{(v)}_{i,j}=\rho^{(v)}_{\tilde{i},\tilde{j}}$ implies that $\{i,j\}=\{\tilde{i},\tilde{j}\}$. Hence,   
$\ell(i \cdot \rho^{(v)}_{i,j}+v_i)=\ell(j \cdot \rho^{(v)}_{i,j}+v_j) \leq \ell(k_0 \cdot \rho^{(v)}_{i,j}+v_{k_0})$ for all $\ell \in N_{\mathfrak{P_0}}(\rho^{(v)}_{i,j})$. Since $N_{\mathfrak{P_0}}(\rho^{(v)}_{i,j}) \cap C$ is full-dimensional,  by the second continuation principle (Principle \ref{seccont_princ}), we deduce that $\ell(i \cdot \rho^{(v)}_{i,j}+v_i)=\ell(j \cdot \rho^{(v)}_{i,j}+v_j) \leq \ell(k_0 \cdot \rho^{(v)}_{i,j}+v_{k_0})$ for all $\ell \in C$. By the first continuation principle (Principle \ref{firstcont_princ}), we deduce that these inequalities are also satisfied by every element in $N_{\mathfrak{P_0}}(\rho^{(v)}_{i,j}) \odot C$. But, this contradicts the existence of $\ell_0$. Hence, the set of vertices of $\phi_{\rm gen}(\tilde{\mathfrak{P}_0})$ is equal to that of $\phi_{\rm gen}(\mathfrak{P}_0)$. Furthermore, the set of vertices of each summand of $\phi_{\rm gen}(\tilde{\mathfrak{P}_0})$ contains the set of vertices of the corresponding summand of $\phi_{\rm gen}(\mathfrak{P}_0)$. Hence, every vertex of  $\phi_{\rm gen}(\tilde{\mathfrak{P}_0})$ is shared by at least two distinct summands. Thus, we obtain a contradiction to the Minkowski-Weyl minimality of $\mathfrak{P}_0$.  \end{proof}

Suppose that $\mathfrak{f}$ is a facet  of the normal cone of a unique vertex $\gamma=\rho^{(v)}_{i,j}$  of $\mathfrak{P}_0$. By Proposition \ref{minweylstruct_prop}, $\mathfrak{f}$ is exclusively either a facet of $N_M(v)$ or a facet of a  maximal cone $C \subseteq \mathcal{F}_v$ that is not contained in the support of $\mathfrak{P}_0$. Suppose that $\mathfrak{f}$ satisfies the latter property, then we have the following proposition.


\begin{proposition}\label{labelproper_prop}
 The maximal cone $C$ does not have $(i,j)$ as a primary label  or there is a vertex $\rho^{(v)}_{\tilde{i},\tilde{j}}$ of $\mathfrak{P}_0$ such that $(i,j,\tilde{i},\tilde{j})$ is not a secondary label of $C$.  
\end{proposition}

\begin{proof}
Suppose for the sake of contradiction that $C$ contains the primary label $(i,j)$ and that $(i,j,\tilde{i},\tilde{j})$ is a secondary label of $C$ for every vertex $\rho^{(v)}_{\tilde{i},\tilde{j}}$ of $\mathfrak{P}_0$.  Consider the collection $\tilde{\mathfrak{P}}_0:=\mathfrak{P}_0 \setminus \{(\gamma,N_{\mathfrak{P_0}}(\gamma))\} \cup \{( \gamma, N_{\mathfrak{P_0}}(\gamma) \odot C)\}$, i.e. the collection obtained from $\mathfrak{P}_0$ by replacing the normal cone of $\gamma$ by its Minkowski sum  with $C$. We start by noting that, by the definition of secondary label, $\sigma(\gamma) \leq \sigma(\rho^{(v)}_{(\tilde{i},\tilde{j})})$ for every $\sigma$ in $N_{\mathfrak{P_0}}(\gamma) \cup C$ and  for every vertex $\rho^{(v)}_{\tilde{i},\tilde{j}}$ of $\mathfrak{P}_0$. By the first continuation principle, this inequality holds for every $\sigma \in N_{\mathfrak{P_0}}(\gamma) \odot C$. The inequalities, required for the vertex-cone collection property of $\tilde{\mathfrak{P}}_0$, corresponding to the other vertices of  $\tilde{\mathfrak{P}}_0$  are inherited from $\mathfrak{P}_0$.  Thus, we conclude that $\tilde{\mathfrak{P}_0}$ is a vertex-cone collection.

We claim that  $\tilde{\mathfrak{P}_0}$ is a solution to $\phi_{\rm gen}$. We first show that the vertex sets of  $\phi_{\rm gen}(\mathfrak{P_0})$  and  $\phi_{\rm gen}(\tilde{\mathfrak{P}_0})$ are equal.  By the definition of Minkowski sum of vertex-cone collections, the vertex set of the $k$-th summand of $\phi_{\rm gen}(\mathfrak{P_0})$ is contained in the vertex set of the $k$-th summand of  $\phi_{\rm gen}(\tilde{\mathfrak{P}_0})$ for all $k \in {\rm Supp}(\phi_{\rm gen})$.  This along with the definition of convex hull (of vertex-cone collections) implies that the vertex set of $\phi_{\rm gen}(\mathfrak{P}_0)$ is contained in the vertex set of $\phi_{\rm gen}(\tilde{\mathfrak{P}_0})$. Furthermore, any vertex of $\phi_{\rm gen}(\tilde{\mathfrak{P}_0})$ that is not a vertex of $\phi_{\rm gen}(\mathfrak{P}_0)$ is of the form $v_{i_l}+i_{l} \cdot \gamma$ and  there is a full-dimensional subcone $\tilde{C}$ of $N_{\mathfrak{P_0}}(\gamma) \odot C$ such that $\sigma(v_{i_l}+i_l \cdot \gamma) \leq \sigma(v_{i_o}+i_o \cdot \gamma)$ for all $o$ from $1$  to $r+1$  and all $\sigma \in \tilde{C}$.  On the other hand, since $C$ has $(i,j)$ as a primary label, we have $\ell(v_i+i \cdot \gamma) \leq \ell(v_{i_o}+i_o \cdot \gamma)$ for all $\ell \in C$ and for all $o$ from $1$ to $r+1$. Furthermore, since $N_{\mathfrak{P_0}}(\gamma) \subseteq \mathcal{F}_v$ and $\phi_{\rm gen}$ is generic, we have $\ell(v_i+i \cdot \gamma) \leq \ell(v_{i_o}+i_o \cdot \gamma)$ for all $o$ from $1$  to $r+1$ and all $\ell \in N_{\mathfrak{P_0}}(\gamma)$. By the first continuity principle, these inequalities are valid over for every $\ell \in N_{\mathfrak{P_0}}(\gamma) \odot C$.  Hence, we deduce that  $v_{i_l}+i_l \cdot \gamma=v_i+i \cdot \gamma=v_j+j \cdot \gamma$. 
We conclude that the vertex sets of  $\phi_{\rm gen}(\mathfrak{P_0})$  and  $\phi_{\rm gen}(\tilde{\mathfrak{P}_0})$ are equal and that $\tilde{\mathfrak{P}_0}$ is a solution to $\phi_{\rm gen}$.
This contradicts the Minkowski-Weyl minimality of $\mathfrak{P}_0$.
\end{proof}

\begin{remark}\rm
In the proofs of Propositions \ref{minweylstruct_prop} and \ref{labelproper_prop}, the generic-local setting ensures that the solution property is preserved in the passage from 
$\mathfrak{P}_0$ to $\tilde{\mathfrak{P}}_0$ via the Minkowski sum operation. In other words, this setting allows the application of the first continuation principle. 
\qed
\end{remark}



Using Proposition \ref{minweylstruct_prop}, we associate a local compatible system to $\mathfrak{P}_0$ as follows. For every vertex $\gamma$ of $\mathfrak{P}_0$, let $C_{1,\gamma},\dots,C_{k_{\gamma},\gamma}$ be the collection of maximal cones of $\mathcal{F}_v$ that are contained in $N_{\mathfrak{P}_0}(\gamma)$. 
Since $\mathfrak{P}_0$ is a solution to $\phi_{\rm gen}$ and $\phi_{\rm gen}$ is generic,  there exist exactly two distinct indices $i_\gamma,~j_\gamma$, say in ${\rm Supp}(\phi_{\rm gen})$ such that for every $\sigma \in \cup_{l} C_{l,\gamma}$, 
 \begin{equation}\label{sigma_eq} \sigma(v_{i_\gamma}+i_{\gamma} \cdot \gamma)=\sigma(v_{j_\gamma}+j_{\gamma} \cdot \gamma) \leq  \sigma(v_{k}+k \cdot \gamma) \end{equation}  
for every $k \in {\rm Supp}(\phi_{\rm gen})$. Consider the collection $\mathcal{C}=\{C_{l,\gamma}\}_{\gamma,l}$ where $\gamma$ varies over all vertices of $\mathfrak{P}_0$ and  $l$ varies over $[1,\dots,k_{\gamma}]$ along with the corresponding set of tuples $\mathcal{I}=\{(i_{\gamma},~j_{\gamma})\}_{\gamma,l}$ (each cone $C_{l,\gamma}$ has the same tuple $(i_{\gamma},~j_{\gamma})$).  In the following proposition, we show that $(\mathcal{C},\mathcal{I})$ is a local compatible system.  

\begin{proposition}\label{vcctolcs_prop} The pair $(\mathcal{C},\mathcal{I})$ is a local compatible system. \end{proposition}

\begin{proof} We verify each property in the definition of local compatible system. Property \ref{comp0_con} is satisfied since the normal cone of every vertex of  $\mathfrak{P}_0$ is a convex polyhedral cone. Property \ref{comp1_con} follows from Equation (\ref{sigma_eq}).  Property \ref{comp2_con} is a consequence of the vertex-cone collection inequalities of $\mathfrak{P}_0$.
Property \ref{comp4_con} follows from Proposition \ref{labelproper_prop}. 
Hence, we conclude that $(\mathcal{C},\mathcal{I})$ is a local compatible system. \end{proof}



The remaining part of the proof of Theorem \ref{genminwey_theo} is to observe that the vertex-cone collection associated to $(\mathcal{C},\mathcal{I})$ is $\mathfrak{P}_0$, i.e the correspondence from a local compatible system to a vertex-cone collection and its opposite correspondence are inverses to each other. This is relatively straightforward and completes the proof of Theorem \ref{genminwey_theo}.  As a corollary, we deduce the following finiteness result.  Fix a generic polyhedral polynomial $\phi_{\rm gen}$ in $\mathcal{A}[Y]$ and a vertex $v$ of the Minkowski sum of its coefficients.  A \emph{$v$-local} solution to $\phi_{\rm gen}$ is a (VCC) solution supported in the normal cone of $v$.

\begin{corollary}
For any vertex $v$ of $M$, the generic polyhedral polynomial $\phi_{\rm gen}$ has finitely many Minkowski-Weyl minimal vertex-cone collection $v$-local solutions.
\end{corollary}

\subsection{Completions, Associated Polyhedra and Shrinking}\label{comasspshr_subsect}

Theorem \ref{minweychar_theointro} asserts that, in the generic-local setting, every Minkowski-Weyl minimal polyhedral solution can be realised in terms of a local compatible system.  On the other hand, by Theorem \ref{genminwey_theo}, any local compatible system corresponds to a Minkowski-Weyl minimal VCC solution.  Keeping this in mind, a key ingredient in the proof of Theorem \ref{minweychar_theointro} is the association of a Minkowski-Weyl minimal VCC solution to any polyhedral solution. We treat this in the following.

Let $\phi_{\rm gen}$ be a generic polyhedral polynomial and let $v$ be a vertex of $M$ (the Minkowski sum of the coefficients of $\phi_{\rm gen}$). Let $P_0$ be a polyhedral solution to $\phi_{\rm gen}$ such that its normal cone  $N(P_0) \subseteq N_M(v)$.

\begin{proposition}\label{normmax_prop}
The normal cone of any vertex of $P_0$ is a union of subcones of  maximal cones of $\mathcal{F}_v$. Furthermore, for any maximal cone $C$ of $\mathcal{F}_v$, the inequality $|\{\gamma \in V(P_0)|~N_{P_0}(\gamma) \cap C \text{~ is full-dimensional} \}| \leq 1$ holds.
\end{proposition}

The proof of Proposition \ref{normmax_prop} is straightforward and is hence, omitted. 
Proposition \ref{normmax_prop} allows us to define the notion of completion of $P_0$ as follows.

\begin{definition}{\rm {(\bf Completion})}
The completion ${\rm Com}(P_0)$ of $P_0$ is the collection $\{ (\gamma, \bar{C}_{\gamma}) \}_{\gamma \in V(P_0)}$ where $\bar{C}_{\gamma}:=\odot_{i=1}^{k} C_{i,\gamma}$ and $C_{1,\gamma},\dots,C_{k,\gamma}$ are the maximal cones of $\mathcal{F}_v$ whose intersection with the  inner normal cone of $\gamma$ is full-dimensional.
\end{definition}

\begin{proposition}\label{comp_prop}
The completion ${\rm Com}(P_0)$ of $P_0$ is a vertex-cone collection and is a solution to $\phi_{\rm gen}$. 
\end{proposition}
\begin{proof}
({\bf VCC property}) Since $P_0$ is a solution to $\phi_{\rm gen}$, by Proposition \ref{minweylstruct_prop} any vertex of $P_0$ is of the form $\rho_{i,j}^{(v)}$ for suitable integers $i \neq j$. 
Let $\gamma_1$ be such a vertex.   The inequality $\sigma(\gamma_1) \leq \sigma(\gamma_2)$ holds for all $\sigma \in N_{P_0}(\gamma_1)$ and for all vertices $\gamma_2 \neq \gamma_1$  of $P_0$.   Since $C_{l,\gamma_1} \cap N_{P_0}(\gamma_1)$ is full-dimensional for all valid $l$ and $C_{l,\gamma_1}$ is a maximal subcone of $\mathcal{F}_v$, the second continuation principal  implies that  $\sigma(\gamma_1) \leq \sigma(\gamma_2)$ for all $\sigma \in C_{l,\gamma_1}$.  By the first continuation principle, we conclude that $\sigma(\gamma_1) \leq \sigma(\gamma_2)$ for all $\sigma \in \bar{C}_{\gamma_1}$. Hence, we conclude that ${\rm Com}(P_0)$ is a VCC.

({\bf Solution Property}) This proof follows the same pattern as the corresponding proofs in Propositions \ref{minweylstruct_prop} and \ref{labelproper_prop}. Hence, we only outline the key steps. 
First, we show that $V(\phi_{\rm gen}({\rm Com}(P_0)))=V(\phi_{\rm gen}(P_0))$. This involves using the genericity of $\phi_{\rm gen}$ along with the two continuation principles.  
The next step is to show for any $i \in {\rm Supp}(\phi_{\rm gen})$, the set of vertices of the $i$-th summand of $\phi_{\rm gen}({\rm Comp}(P_0))$ contains the set of vertices of the $i$-th summand of  $\phi_{\rm gen}(P_0)$. Since $P_0$ is a solution to $\phi_{\rm gen}$, we deduce that ${\rm Com}(P_0)$ is also a solution to $\phi_{\rm gen}$.
 \end{proof}
 
 Note that the notion of completion can also be analogously defined for vertex-cone collection solutions to $\phi_{\rm gen}$ whose support is contained in $N_M(v)$.



\begin{proposition}\label{minvcc_prop}
Fix a vertex $v$ of $M$. For any $v$-local (vertex-cone collection) solution $\mathfrak{B}_0$ to $\phi_{\rm gen}$, there exists a $v$-local Minkowski-Weyl minimal (vertex-cone collection) solution with the same vertex set as $\mathfrak{B}_0$.
\end{proposition}

\begin{proof}
Consider the set $\{\tilde{\mathfrak{B}_0}\}$ of all $v$-local solutions to $\phi_{\rm gen}$ whose support contains that of $\mathfrak{B}_0$.  By Proposition \ref{comp_prop} and the finiteness of the number of maximal cones of $\mathcal{F}_v$, the set $\{{\rm Com}(\tilde{\mathfrak{B}_0})\}$ is finite. Any element in this set with maximal support is a $v$-local Minkowski-Weyl minimal solution with the same vertex set as $\mathfrak{B}_0$.
\end{proof}

{\bf Associated Polyhedra of an LCS:} In the following, we define another ingredient of Theorem \ref{minweychar_theointro}: the associated polyhedron of the restriction of an LCS (to a convex polyhedral subcone of its support).  Let $(\mathcal{C},\mathcal{I})$ be a local compatible system and let $L$ be a convex polyhedral subcone of its support. Suppose that $\tilde{\mathfrak{T}}$ is the subcone of $\mathfrak{T}$ consisting of elements $\gamma$ such that $C_{\gamma} \cap L$ is full-dimensional.  By definition of an LCS, there is a unique polyhedron with vertex set $\tilde{\mathfrak{T}}$ and with the vertex $\gamma$ having normal cone $C_{\gamma} \cap L$ for each $\gamma \in \tilde{\mathfrak{T}}$. We refer to this polyhedron as the \emph{associated polyhedron of the restriction of $(\mathcal{C},\mathcal{I})$ to $L$}.

{\bf A Shrinking Principle:} The following shrinking principle will be invoked in the proof of Theorem \ref{minweychar_theointro}.

\begin{principle}{\rm ({\bf Shrinking Principle})}\label{shrin_prin}
If a vertex-cone collection $\mathfrak{B}_0$ is a solution to a polyhedral polynomial $\phi$, then so is any vertex-cone collection whose support is contained in that of $\mathfrak{B}_0$. 
\end{principle}

Principle \ref{shrin_prin} follows from the underlying definitions.   For a polyhedron $P$, let ${\rm Supp}(P)$ denote its inner normal cone, i.e. the union of inner normal cones of all its vertices.

\subsection{Proof of Theorem \ref{minweychar_theointro}}

($\Rightarrow$) Consider the completion ${\rm Com}(P_0)$ of $P_0$.  By Proposition \ref{comp_prop},  ${\rm Com}(P_0)$ is a VCC solution to $\phi_{\rm gen}$.  Let ${\rm Com}^{\rm min}(P_0)$ be the Minkowski-Weyl minimal VCC solution associated with ${\rm Com}(P_0)$ (Proposition \ref{minvcc_prop}).  By Theorem \ref{genminwey_theo}, there is a local compatible system $(\mathcal{C}_0,\mathcal{I}_0)$ associated with ${\rm Com}^{\rm min}(P_0)$.  We claim that $P_0$ is the associated polyhedron  of  $(\mathcal{C}_0,\mathcal{I}_0)$  restricted to a maximal polyhedral convex subcone of its support (with ``support'' as defined in Remark \ref{supplcs_rem}). Suppose not, then there is a polyhedral convex subcone of the support of $(\mathcal{C}_0,\mathcal{I}_0)$  that contains $N(P_0)$. Let $\tilde{P}_0$ be its associated polyhedron. Note that the vertex sets of $\tilde{P}_0$ and $P_0$ are equal, and $\tilde{P}_0$ is a polyhedral solution to $\phi_{\rm gen}$ (Principle \ref{shrin_prin}). This contradicts the Minkowski-Weyl minimality of $P_0$. 

($\Leftarrow$) The associated polyhedron $P_0$ of the restriction of a local compatible system  $(\mathcal{C}_0,\mathcal{I}_0)$ (at $v$) to a maximal polyhedral convex set is a solution to $\phi_{\rm gen}$ (Principle \ref{shrin_prin}). Furthermore, $P_0$ must be a Minkowski-Weyl minimal polyhedral solution. Assume not for the sake of contradiction, then there is another polyhedral solution $\tilde{P}_0$ such that $N(P_0) \subset N(\tilde{P_0})$ and $V(P_0)=V(\tilde{P_0})$. Comparing their completions ${\rm Com}(P_0)$ and ${\rm Com}(\tilde{P_0})$, we note that ${\rm Supp}({\rm Com}(P_0)) \subseteq {\rm Supp}({\rm Com}(\tilde{P_0}))$. If ${\rm Supp}({\rm Com}(P_0))={\rm Supp}({\rm Com}(\tilde{P_0}))$, then ${\rm Com}(P_0)={\rm Com}(\tilde{P_0})$. By construction, the vertex-cone collection  associated with $(\mathcal{C}_0,\mathcal{I}_0)$ (via Theorem \ref{genminwey_theo}) is ${\rm Com}(P_0)$ and  this contradicts the fact that $N(P_0)$ is a maximal convex subcone of ${\rm Supp}({\rm Com}(P_0))$. Hence, ${\rm Supp}({\rm Com}(P_0)) \subset {\rm Supp}({\rm Com}(\tilde{P_0}))$. But then, $V(P_0)=V(\tilde{P_0})=V({\rm Com}(P_0))=V({\rm Com}(\tilde{P_0}))$. 
But, this contradicts the Minkowski-Weyl minimality of ${\rm Com}(P_0)$. Hence, we conclude that $P_0$ is a Minkowski-Weyl minimal polyhedral solution to $\phi_{\rm gen}$.

\section{A Local-Global Principle}\label{locglo_sect}

Having studied local solutions to a polyhedral polynomial, we move to the global setting. Recall that a typical polyhedral polynomial is given by $\phi(Y)=Q_{d} \odot Y^{d} \oplus Q_{i_r} \odot  Y^{i_r} \oplus \dots \oplus Q_{i_0} \odot Y^{i_0}$ and $M:=Q_d \odot Q_{i_r} \odot \cdots \odot Q_{i_0}$ is the Minkowski sum of its coefficients. 

\begin{definition}{\rm {(\bf Complete Local and Global Solutions})}
For a vertex $v$ of $M$, a $v$-local solution to $\phi$ is called complete if it has maximum support, i.e. support equal to $N_M(v)$.
A polyhedral solution $P_0$ to $\phi$ is called a global solution 
if the normal fan of $P_0$ contains the normal fan of $M$.
\end{definition}

If $P_0$ is a global solution, then for any vertex $v$ of $M$, the polyhedron $P_0 \odot (N_M(v))^{\star}$ 
is a $v$-local solution. The central question in this context is the following converse: 

Suppose that a polyhedral polynomial $\phi$ has a complete $v$-local solution for every vertex $v \in M$, under what circumstances does it have a global solution?

The following simple examples shed light on this question.

\begin{example}\label{locglonoglo_ex}
\rm 
Consider the linear polynomial $\phi=[-1,1] \odot Y \oplus 0$. In this case, $M=[-1,1]$, and the inner normal fans $N_M(-1)=\mathbb{R}_{\geq 0}$ and $N_M(1)=\mathbb{R}_{\leq 0}$. 
Complete local solutions exist for both the vertices of $M$: the affine cone $1 \odot \mathbb{R}_{\geq 0}$ is a complete $(-1)$-local solution and $-1 \odot \mathbb{R}_{\leq 0}$ is a complete $1$-local solution. However, comparing the volumes of the two coefficients tells us that there no global solutions. \qed
\end{example}

\begin{example}\label{locgloglo_ex}
\rm Consider the linear polynomial $\phi=[-1,1] \odot Y \oplus [-2,2]$. Here, $M=[-3,3]$, and the normal fans $N_M(-3)=\mathbb{R}_{\geq 0}$ and $N_M(3)=\mathbb{R}_{\leq 0}$.
Complete local solutions exist for both $-3$ and $3$: the affine cone $-1 \odot \mathbb{R}_{\geq 0}$ is a complete $(-3)$-local solution and $1 \odot \mathbb{R}_{\leq 0}$ is a complete $3$-local solution.  Furthermore, $[-1,1]$ is a global solution.
\qed
\end{example}

The key difference between Examples \ref{locglonoglo_ex} and \ref{locgloglo_ex} is that the local solutions are ``oriented'' so that they cannot be ``glued" to form a polytope in the former and can be glued to form a polytope in the latter.  This suggests an additional orientation condition that allows for such a gluing.  This leads to the Property (P) in Theorem \ref{locglo_theo}. 


{\bf Proof of Theorem \ref{locglo_theo}:} Analogous with the classical local-global principle,  the passage from global to local $(\Rightarrow)$ is more straightforward.

 ($\Rightarrow$)  Suppose that $\hat{P_0}$ is a global solution to $\phi$. For every vertex $v$ of $M$, the polyhedron $S_v=\hat{P_0} \odot (N_M(v))^{\star}$ is a complete $v$-local solution.  By construction, every vertex $\gamma$ of $S_v$ is a vertex of $\hat{P_0}$ and $N_{S_v}(\gamma)=N_{\hat{P}_0}(\gamma) \cap N_M(v)$. As a consequence,  the Property (P) is satisfied.
 
 ($\Leftarrow$)  Suppose that $w$ is a vertex of $\phi(P_0)$. A generic element in $N_{\phi(P_0)}(w)$ attains unique minima over both $M$ and $P_0$. Suppose that $\ell$ is such an element, and  let $v$ and $\gamma$ be the unique vertices in $M$ and $P_0$, respectively where $\ell$ attains its minimum.  Furthermore,  $N_M(v) \cap N_{P_0}(\gamma)$ is full-dimensional. 
 Hence, by Property \ref{ori_prop},  $\gamma$ is a vertex of $S_v$.  Since $V(S_v) \subseteq V(P_0)$,  we deduce that $\ell$ attains a unique minimum over $S_v$ at $\gamma$.
 Consider $\phi(S_v)$. Note that $\ell \in N(\phi(S_v))$ and that it attains a unique minimum over $\phi(S_v)$ at $w$. We deduce that $w$ is a vertex of $\phi(S_v)$. Hence, $w$ is shared by at least two distinct summands of $\phi(S_v)$. Suppose that $i$ and $j$ are two such summands.  We have $w=v_i+i \cdot \gamma=v_j+j \cdot \gamma$ (recall that $v=v_{i_0}+v_{i_1}+\dots+v_{i_{r+1}}$ is the Minkowski decomposition of $v$). Since $\ell$ attains its unique minimum at $v_i,~v_j$ and $\gamma$ over  $Q_i,~Q_j$ and $P_0$, respectively, we deduce that $w$ is a vertex of the $i$-th and the $j$-th summand of $\phi(P_0)$. This shows that $P_0$ is a solution to $\phi$.  \qed

\begin{remark} 
\rm
Example \ref{locgloglo_ex} satisfies Property \ref{ori_prop} but not Example \ref{locglonoglo_ex}. \qed
\end{remark}

Let $P_0={\rm CH}(v_0-v_1|~v_0+v_1 \in V(M))$ where $v_0+v_1$ is the Minkowski decomposition.   Specialising to the case where both the coefficients are polytopes, Theorem \ref{locglo_theo} yields the following characterisation of Minkowski summands. 

\begin{proposition}\label{minsum_prop}
The polytope $Q_1$ is a Minkowski summand of the polytope $Q_0$ if and only if every point of the form $v_0-v_1$ is a vertex of $P_0$ and for every $\gamma \in V(P_0)$, its normal cone $N_{P_0}(\gamma)=\cup_{v|~\rho^{(v)}_{0,1}=\gamma}N_M(v)$ where $v \in V(M)$.
\end{proposition}

In the following, we use Proposition \ref{minsum_prop} to derive a criterion for weak Minkowski summands, due to Shephard, see \cite{CaDoGoRoYi22} for a recent work on this topic. This criterion is foundational to the theory of Minkowski weights \cite{FulStu97}.  
 Recall that a polytope $Q_1$ is called a \emph{weak Minkowski summand} of a polytope $Q_0$ if there is a positive real $\lambda_0$ such that $Q_1$ is a Minkowski summand of  $\lambda_0 \cdot Q_0$. Shephard's criterion states that $Q_1$ is a weak Minkowski summand of $Q_0$ if and only if there is a map ${\rm Sp}: V(Q_0) \rightarrow V(Q_1)$ such that for every edge $(u_0,v_0)$  of $Q_0$, there is a non-negative real $\lambda$ such that $\lambda \cdot (u_0-v_0)=({\rm Sp}(u_0)-{\rm Sp}(v_0))$.

\subsection{Shephard's criterion  via Proposition \ref{minsum_prop}} 
If $Q_1$ is a weak Minkowski summand of $Q_0$ with $\lambda_0$ as a corresponding dilation, then each vertex $\lambda_0 \cdot v_0$ of $\lambda_0 \cdot Q_0$ has a unique Minkowski decomposition of the form $v_1+\hat{v}$ where $v_1$ is a vertex of $Q_1$. Define the map ${\rm Sp}$ as ${\rm Sp}(v_0):=v_1$.  By Proposition \ref{minsum_prop}, we have the following: 
\begin{enumerate}
\item\label{locglo_con1} Every point $\lambda_0 \cdot v_0-v_1$ is a vertex of the polytope $\tilde{P_0}:={\rm CH}(\lambda_0 \cdot v_0-v_1|~\lambda_0 \cdot v_0+v_1 \in V(\tilde{M}))$ where $\tilde{M}:=(\lambda_0 \cdot Q_0) \odot Q_1$.

\item\label{locglo_con2} For every vertex $\gamma$ of $\tilde{P_0}$, its normal cone $N_{\tilde{P_0}}(\gamma)=\cup_{v|~\rho^{(v)}_{0,1}=\gamma}N_{\tilde{M}}(v)$ where $v \in V(\tilde{M})$.
 \end{enumerate}
 
 Consider an edge $(u_0,v_0)$ of $Q_1$. Suppose that ${\rm Sp}(u_0)={\rm Sp}(v_0)$, then $\lambda_{u_0,v_0}=0$ satisfies the required condition. On the other hand, if ${\rm Sp}(u_0) \neq {\rm Sp}(v_0)$ but $\lambda_0 \cdot u_0- {\rm Sp}(u_0)=\lambda_0 \cdot v_0-{\rm Sp}(v_0)$, then set $\lambda_{u_0,v_0}$ to $\lambda_0$. 
 Suppose that ${\rm Sp}(u_0) \neq {\rm Sp}(v_0)$ and 
$\lambda_0 \cdot u_0- {\rm Sp}(u_0) \neq \lambda_0 \cdot v_0-{\rm Sp}(v_0)$. Since the normal fans of $Q_0$ and $\tilde{M}$ are the same, $(\lambda_0 \cdot u_0+{\rm Sp}(u_0), \lambda_0 \cdot v_0+{\rm Sp}(v_0))$ is an edge of $\tilde{M}$, call it $e^+$. By Property \ref{locglo_con1},  $\lambda_0 \cdot u_0-{\rm Sp}(u_0)$ and $\lambda_0 \cdot v_0+{\rm Sp}(v_0)$ are both vertices of $\tilde{P_0}$. By Property \ref{locglo_con2}, $N_{\tilde{P_0}}(\lambda_0 \cdot u_0-{\rm Sp}(u_0)) \supseteq N_{\tilde{M}}(\lambda_0 \cdot u_0+{\rm Sp}(u_0))$ and $N_{\tilde{P_0}}( \lambda_0 \cdot v_0-{\rm Sp}(v_0)) \supseteq N_{\tilde{M}}( \lambda_0 \cdot v_0+{\rm Sp}(v_0))$. Hence, $(\lambda_0 \cdot u_0-{\rm Sp}(u_0), \lambda_0 \cdot v_0-{\rm Sp}(v_0))$ is an edge of $\tilde{P}_0$, call it $e^{-}$.  Furthermore, the normal cone of $e^{+}$ is contained in that of $e^{-}$. Thus, the hyperplanes spanned by these normal cones coincide. We conclude that there is a 
$\bar{\lambda} \in \mathbb{R}$ such that $\bar{\lambda} \cdot ((\lambda_0 \cdot u_0-{\rm Sp}(u_0))- (\lambda_0 \cdot v_0-{\rm Sp}(v_0)))=(\lambda_0 \cdot u_0+{\rm Sp}(u_0))-(\lambda_0 \cdot v_0+{\rm Sp}(v_0))$. Furthermore, since $N_{\tilde{P_0}}(\lambda_0 \cdot u_0-{\rm Sp}(u_0)) \supseteq N_{\tilde{M}}(\lambda_0 \cdot u_0+{\rm Sp}(u_0))$, we deduce that $\bar{\lambda} \geq 0$.
Rearranging terms yields $\hat{\lambda} \cdot (u_0-v_0)=({\rm Sp}(u_0)-{\rm Sp}(v_0))$ where $\hat{\lambda}=\lambda_0 \frac{\bar{\lambda}-1}{\bar{\lambda}+1}$.
We finally note that $\hat{\lambda} \geq 0$ since the normal cones of $u_0$ and ${\rm Sp}(u_0)$ intersect in full dimension. 


 
 


\section{Examples}\label{example_sect}

We  illustrate the theory developed so far for polyhedral polynomials of small degree (up to cubics). 


\subsection{Linear Case}
Consider a linear polynomial of the form $Q_1 \odot Y \oplus Q_0$ where $Q_1,Q_0 \in \mathcal{A}$,   also treated in \cite[Section 6]{Val22}.  Let $v$ be a vertex of $M(:=Q_1 \odot Q_0)$  and let  $v_1+v_0$ be its Minkowski decomposition.  Since $\phi$ has only two summands, a complete $v$-local solution exists for every $v$ and it is $(v_0-v_1) \odot (N_M(v))^{\star}$. The associated local compatible system is $\{N_M(v)\},\{(0,1)\}$.

In the light of Remark \ref{minsum_rem}, global solutions need not exist. Consider the case where $Q_1$ and $Q_0$ are polytopes such that their volumes ${\rm vol}(Q_1)>{\rm vol}(Q_0)$. Any polytope $P_0$ that is a solution to $\phi$ must satisfy  $Q_1 \odot P_0=Q_0$. On the other hand, ${\rm vol}(Q_1 \odot P_0) \geq {\rm vol}(Q_1)$. This obstructs the existence of polytopal solutions.  We provide another example to illustrate Theorem  \ref{locglo_theo}.  Consider the inverted variant of the icosahedral example from the introduction, i.e. $Q_1 \odot Y \oplus Q_0$ where $Q_0$ is that regular icosahedron and $Q_1$ is the edge $(v_1,v_2)$. We claim that this equation does not have a global solution. To prove this via Theorem \ref{locglo_theo}, suppose that there exists a polytopal solution $P_0$. The origin $O$ must be a vertex of $P_0$ since it is $\rho^{(2 v_1)}_{0,1}=\rho^{(2 v_2)}_{0,1}$. Consider $N_{P_0}(O)$. Since the union of the normal cones of $2v_1$ and $2v_2$ is not convex, there exists a vertex $v \notin \{2v_1,2v_2\}$ of $M$ such that $N_{M}(v) \cap N_{P_0}(O)$ is full-dimensional. But, $\rho_{0,1}^{(v)} \neq O$ and thus, the criterion for existence of a global solution in Theorem \ref{locglo_theo} is not satisfied.

 \subsection{Quadratics and Cubics}
 

{\bf Quadratics:} Consider a polynomial $\phi$ of the form $Q_2 \odot Y^2 \oplus Q_1 \odot Y \oplus Q_0$ where $Q_2,Q_1,Q_0 \in \mathcal{A}$. Let $M=Q_2 \odot Q_1 \odot Q_0$. 
 Recall that any vertex $v$ of $M$ has a unique Minkowski decomposition $v=v_2+v_1+v_0$ where $v_k \in V(Q_k)$. Let $\Delta_v:=2v_1-v_0-v_2$.  Let $H_{\rm disc}(v)$ be the subspace $\{\ell \in (\mathbb{R}^n)^\star |~\ell(\Delta_v)=0\}$, and let $H^{+}_{\rm disc}(v):=\{\ell \in (\mathbb{R}^n)^\star |~\ell(\Delta_v) \geq 0\}$ and $H^{-}_{\rm disc}(v):=\{\ell \in (\mathbb{R}^n)^\star |~\ell(\Delta_v) \leq 0\}$. Since $\Delta_v$  plays a role analogous to the discriminant, we refer to it as the \emph{discriminant at $v$} and the corresponding spaces as \emph{discriminantal subspaces} (at $v$).  The following proposition classifies the $v$-local complete solutions of $\phi$ in terms of discriminantal subspaces at $v$. 

\begin{proposition}
For any vertex $v$ of $M$, a complete $v$-local solution to $\phi$ exists and is precisely one of the following types: 
\begin{itemize}
\item {\rm ({\bf  Degenerate})}  The vector $\Delta_v=0$. In this case, $\tau:=\rho^{(v)}_{0,1}=\rho^{(v)}_{1,2}=\rho^{(v)}_{0,2}$ and $\tau \odot (N_M(v))^{\star}$ is the only $v$-local complete solution.

\item If $N_M(v) \subseteq H^{-}_{\rm disc}(v)$, then $\rho^{(v)}_{0,1} \odot (N_M(v))^{\star}$ and $\rho^{(v)}_{1,2}\odot (N_M(v))^{\star}$ are the only $v$-local complete solutions.

\item If  $N_M(v) \subseteq H^{+}_{\rm disc}(v)$, then $\rho^{(v)}_{0,2} \odot  (N_M(v))^{\star}$ is the only $v$-local complete solution.

\item If $N_M(v)$ intersects both $H^{-}_{\rm disc}(v)$ and $H^{+}_{\rm disc}(v)$ in full dimension, then $(\rho^{(v)}_{1,2}\odot (N_M(v) \cap H^{-}_{\rm disc}(v))^{\star}) \cap (\rho^{(v)}_{0,2} \odot  (N_M(v) \cap H^{+}_{\rm disc}(v))^{\star})$ is the only $v$-local complete solution.
\end{itemize}
\end{proposition} 

Gluing local solutions to a global one seems involved. 







{\bf Reduced Cubics:} Consider a cubic $\phi$ of the form $Q_3 \odot Y^3 \oplus Q_1 \odot Y \oplus Q_0$.  Let $v$ be a vertex of $M:=Q_3 \odot Q_1 \odot Q_0$ with Minkowski decomposition $v_3+v_1+v_0$. Analogous to the quadratic case, let $\Delta_v:=3v_1-2v_0-v_3$ be the discriminant of $\phi$ at $v$. The spaces $H_{\rm disc}(v), H^{+}_{\rm disc}(v)$ and $H^{-}_{\rm disc}(v)$ are defined analogously and the $v$-local complete solutions can be characterised as follows. 

\begin{proposition}
For any vertex $v$ of $M$, a complete $v$-local solution to $\phi$ exists and is precisely one of the following types:
\begin{itemize}
\item  {\rm ({\bf  Degenerate})} The vector $\Delta_v=0$.  In this case,  $\tau:=\rho^{(v)}_{0,1}=\rho^{(v)}_{1,3}=\rho^{(v)}_{0,3}$ and $\tau \odot (N_M(v))^{\star}$ is the only $v$-local complete solution.

\item If $N_M(v) \subseteq H^{-}_{\rm disc}(v)$,  then $\rho^{(v)}_{0,1} \odot (N_M(v))^{\star}$ and $\rho^{(v)}_{1,3}\odot (N_M(v))^{\star}$ are the only $v$-local complete solutions.

\item If  $N_M(v) \subseteq H^{+}_{\rm disc}(v)$, then $\rho^{(v)}_{0,3} \odot  (N_M(v))^{\star}$ is the only $v$-local complete solution. 

\item If $N_M(v)$ intersects both $H^{-}_{\rm disc}(v)$ and $H^{+}_{\rm disc}(v)$ in full dimension, then $(\rho^{(v)}_{1,3}\odot (N_M(v) \cap H^{-}_{\rm disc}(v))^{\star}) \cap (\rho^{(v)}_{0,3} \odot  (N_M(v) \cap H^{+}_{\rm disc}(v))^{\star})$ is the only $v$-local complete solution.

\end{itemize}

\end{proposition}

{\bf A Non-Generic Cubic:} The following example sheds light on some difficulties that arise when the genericity condition in Theorem \ref{minweychar_theointro} is dropped. 
Let $\phi=Q_3 \odot Y^3 \oplus Q_2 \odot Y^2 \oplus Q_1 \odot Y \oplus Q_0$ where $Q_2,Q_1$ and $Q_0$ are polyhedra in $\mathcal{A}$ such that there exists a vertex of  $Q_0 \odot Q_1 \odot Q_2$ with Minkowski decomposition $v_0+v_1+v_2$ such that  $\{v_0,v_1,v_2\}$ is linearly independent and $Q_3=\{v_3\}$ where $v_3=v_2+v_1-v_0$. For example, $Q_2,Q_1$ and $Q_0$ can be the three edges incident on any given vertex of a three-cube.  The point $\hat{v}=v_0+v_2+v_1+v_3$ is a vertex of $M (:=Q_0 \odot Q_1 \odot Q_2 \odot Q_3)$. The genericity condition is not satisfied at this vertex since $\rho^{(\hat{v})}_{1,3}=\rho^{(\hat{v})}_{0,2}$. Analysing the fan $\mathcal{F}_{\hat{v}}$, we find that it has two maximal cones $C_1$ and $C_2$
with primary labels $(1,3)$ and $(0,2)$, respectively and corresponding secondary labels $(1,3,1,2)$ and $(0,2,1,2)$. The set $C_1 \cup C_2$ is not convex and  this non-convexity causes difficulties in the proof of Theorem \ref{genminwey_theo} where Minkowski sums of maximal cones are taken. 






\section{Multiplicity and Discriminants}\label{multdisc_sect}

 In this section, we take $\mathcal{A}=\mathcal{P}_{\omega}$ for a suitable $\omega$.   We start with the notion of multiplicity in this polyhedral setting taking cue from the classical and tropical situations \cite{DicFeiStu07}. The multiplicity of a root of a tropical polynomial is the ``absolute change" in the slope of the polynomial (regarded as a function) \cite{Mor19}. This does not seem to directly generalise to the polyhedral setting since the notion of ``slope of a polyhedral polynomial at a point" is, as far as we are aware, not available. Another formulation of multiplicity of a root in the tropical setting is the following: a root $\lambda$ of a tropical polynomial $\phi$ has multiplicity $m-1$  if $m$ is the maximum integer for which there exists a tropical polynomial $\tilde{\phi}$ that as a function is the same as $\phi$ and $\tilde{\phi}(\lambda)$ is attained by at least $m$ terms. In order to extend this definition to the polyhedral case, we regard a polynomial $\phi \in \mathcal{A}[Y]$ as a function $F_\phi:\mathcal{A} \rightarrow \mathcal{A}$ given by $F_\phi(P)=\phi(P)$. Polynomials $\phi_1,\phi_2 \in \mathcal{A}[Y]$ are said to \emph{represent the same function} if $F_{\phi_1}=F_{\phi_2}$. 
 For a polyhedral polynomial $\phi= \oplus_i (Q_i \odot Y^{i})$ and $\ell \in N(M)$ (here, $M:=\odot_i Q_i$ and N(.) is the inner normal cone), the tropical polynomial $\phi_{\ell}$ is given by $\oplus_i (m_i \odot y^{i})$ where $m_i={\rm min}_{{\bf r} \in Q_i} \ell ({\bf r})$.  In the following,  by a \emph{root} of $\phi \in \mathcal{A}[Y]$, we mean an $\omega$-positive polyhedral root of $\phi$. 
 For simplicity, we do not handle the case where $0_{\mathcal{A}}$ is a root.
 
 
 \begin{proposition}\label{rootinv_prop}
 Suppose that the polynomials $\phi_1,\phi_2 \in \mathcal{A}[Y]$ represent the same function.  Any root $P_0$ of $\phi_1$ is also a root of $\phi_2$.
 \end{proposition}
 \begin{proof}
By evaluating both polynomials at the origin, we note that the Minkowski-Weyl cones of $M_1$ and $M_2$ are the same (as usual $M_i$ is the Minkowski sum of the coefficients of $\phi_i$). Consider $L=M_1 \odot M_2 \odot P_0$. Let $v$ be a vertex of $L$ 
and let $\ell \in N_L(v)$. The tropical polynomials $\phi_{1,\ell}$ and $\phi_{2,\ell}$ are equal as functions and hence, have the same set of roots.  Thus, ${\rm min}_{{\bf r} \in P_0} \ell ({\bf r})$ is a root of $\phi_{2,\ell}$ for all $\ell \in N_L(v)$. Since $N_L(v)$ is full-dimensional, we deduce that $P_0 \odot (N_L(v))^{\star}$ is a root of $\phi_2$ for every vertex of $L$. Hence, $P_0 \odot L^{\rm cone}$ is a root of $\phi_2$ (here, $L^{\rm cone}$ is the Minkowski-Weyl cone of $L$).   The equality of the Minkowski-Weyl cones of $M_1$ and $M_2$ implies that $P_0$ is also a root of $\phi_2$. 
\end{proof}
 
 Proposition \ref{rootinv_prop} justifies the following definition of multiplicity.

\begin{definition}{\rm ({\bf Polyhedral Multiplicity})}
 Let $\phi \in \mathcal{A}[Y]$ and let $P_0$ be a root of $\phi$. The root $P_0$ has multiplicity $m-1$ if $m$ is the maximum integer for which there exists a polynomial $\tilde{\phi} \in \mathcal{A}[Y]$ that represents the same function as $\phi$ and every vertex of $\tilde{\phi}(P_0)$ is shared by at least $m$ summands.
\end{definition}

Classically, the theory of discriminants \cite[Introduction]{GKZ09} studies algebraic conditions on the coefficients of a polynomial whenever it has a root of higher multiplicity, or more generally a ``degeneracy".  As we have already seen while constructing roots, a root of a polyhedral polynomial obscures information about the coefficients of the polynomial unless its support is ``sufficiently large''.  This observation leads to the following notion.


\begin{definition}{\rm ({\bf Degenerate Polynomial)}}
A polynomial $\phi \in \mathcal{A}[Y]$ given by $\oplus_{i \in {\rm Supp}(\phi)} (Q_{i} \odot Y^{i})$ is called weakly degenerate if it has a root $P_0$ of multiplicity at least two and is called degenerate if in addition $P_0$ has support equal to that of $\odot_{i  \in {\rm Supp}(\phi)} Q_{i}$.
\end{definition}

A natural question in this context is whether there is a notion of a ``discriminant''.  For this, factorisation of a polynomial into a product of linear ones is relevant. The classical and tropical situations are as follows.

A non-zero polynomial $p=\sum_{i=0}^{d} c_i x^i \in \mathbb{K}[x]$ where $\mathbb{K}$ is an algebraically closed field can be expressed as $c_d \prod_{i=1}^{d} (x-\alpha_i)$ where $\alpha_1,\dots,\alpha_d$ are the roots of $p$.  An analogous statement does not hold for tropical polynomials in one variable \cite[Pages 9,10]{MacStu15}.  However, a weaker version holds: every tropical polynomial ${\rm min}_{i} (t_i+i \cdot y)$ of degree $d$ is \emph{as a function} equal to $t_d+\sum_{i=1}^{d}{\rm min}(y,\beta_i)$ where $\beta_1,\dots,\beta_d$ are its roots.
We show the following polyhedral  generalisation of this statement. 

\begin{proposition}\label{prodform_prop}
For every polynomial $\phi \in \mathcal{A}[Y]$ of degree $d$ with leading coefficient $Q_d$, there exist $d$ roots $P_1,\dots,P_d$ of $\phi$ such that $Q_d \odot (\odot_{i=1}^{d} (Y \oplus P_i))$ represents the same function as $\phi$.
\end{proposition}

One distinct feature here is that in general, a polyhedral polynomial has infinitely many roots. Hence, a ``choice'' of $d$ roots is necessary. 
We show Proposition \ref{prodform_prop} by generalising a proof of  Kapranov's Theorem in one variable \cite[Section 3.1]{MacStu15}.  

The strategy of the proof is as follows. Consider a polynomial $p=\sum_{i=0}^{d} c_i x^i \in \mathbb{K}_P[x]$ whose polyhedralisation is $\phi$. The polynomial $p$ admits a factorisation as  $c_d \prod_{i=1}^{d} (x-\alpha_i)$. Set $P_i:={\rm New}(\alpha_i) \odot {\rm MWC}(M)$ where $M=\odot_i {\rm New}(c_i)$ and ${\rm MWC}(.)$ is the Minkowski-Weyl cone. Note that ${\rm New}(\alpha_i)$ is, a priori, not a polyhedron but it turns out that $P_i$ is (Corollary \ref{poly_cor}). The argument then proceeds by showing that the evaluations of $\phi$ and $Q_d \odot (\odot_{i=1}^{d} (Y \oplus P_i))$ agree on a ``dense subset'' of $\mathcal{A}$. The ``continuity'' of the two functions $\phi$  and $Q_d \odot (\odot_{i=1}^{d} (Y \oplus P_i))$ implies that they agree on all of $\mathcal{A}$. In order to make sense of these notions of density and continuity, we consider a certain metric on $\mathcal{A}$. We refer to the appendix for more details.

\begin{definition}{\rm ({\bf Product Form})}\label{prodform_def}
A polyhedral polynomial $\phi \in \mathcal{A}[Y]$ of degree $d$ is said to be in product form if there exist polyhedra $Q,P_1,\dots,P_d \in \mathcal{A}$ such that $\phi=Q \odot (\odot_{i=1}^{d} (Y \oplus P_i))$.\end{definition}

Note that in Definition \ref{prodform_def}, the equality is not just as functions but is in the semiring  $\mathcal{A}[Y]$. The tropical discriminant suggests the following definition.


\begin{definition}{\rm ({\bf Polyhedral Discriminant})}\label{polydisc_def}
Let $\Xi$ be a finite subset of $\mathbb{Z}_{\geq 0}$. A polynomial $\Delta_{\Xi} \in \mathcal{A}[Y_i|~i \in \Xi]$ is called a polyhedral discriminant with respect to $\Xi$ if it satisfies the following property: a tuple $(Q_i)_{i \in \Xi}$ is a root of $\Delta_{\Xi}$ if and only if the polynomial $\oplus_{i \in \Xi} (Q_i \odot Y^{i}) \in \mathcal{A}[Y]$ is  in product form and is degenerate.
\end{definition}

 We recall the notion of $X$-discriminant from \cite[Chapter 1]{GKZ09} . Given a finite subset $\Xi$ of $\mathbb{Z}_{\geq 0}$ with maximum element $d$, let $X$ be the image of the map $\mathbb{P}_{\mathbb{C}}^1 \rightarrow \mathbb{P}_{\mathbb{C}}^{|\Xi|-1}$ given by $[x^i y^{d-i}|~i \in \Xi]$.  Its dual variety $X^{\vee}$ is an irreducible, hypersurface in  $(\mathbb{P}^{|\Xi|-1})^{\star}$.  The $X$-discriminant associated with $\Xi$ is the defining equation of $X^{\vee}$.  In the following, we construct a candidate polyhedral discriminant in terms of the $X$-discriminant.
 Before this, note that the coefficients of the $X$-discriminant are complex numbers and hence, their polyhedralisation is well-defined. 

\begin{definition}{\rm ({\bf Polyhedralised Discriminant})}
Let $\Xi$ be a finite subset of $\mathbb{Z}_{\geq 0}$. The polyhedralised discriminant  $\tilde{\Delta}_{\Xi} \in \mathcal{A}[Y_i|~i \in \Xi]$ is the polyhedralisation of the $X$-discriminant associated with $\Xi$. \end{definition}

We, perhaps naively, expect that the polyhedralised discriminant $\tilde{\Delta}_{\Xi}$ is a polyhedral discriminant for every choice of $\Xi$.  In the following, we show one direction of the required equivalence and a weak converse. Before this, we state a principle from linear algebra that plays a key role in both the proofs. Recall that $W$ is a finite dimensional vector space. 

\begin{principle}\label{linalg_princ}
Let $v_1,v_2 \in W$. If there is a full-dimensional cone $C$ in $W^{\star}$ such that $\ell(v_1)=\ell(v_2)$ for all $\ell \in C$, then $v_1=v_2$.
\end{principle}

\begin{proposition}\label{polydisc_prop}
If $\phi=\oplus_{i \in {\rm Supp}(\phi)} (Q_{i} \odot Y^{i}) \in \mathcal{A}[Y]$ is a degenerate polynomial in product form, then its tuple of coefficients $(Q_i|~i \in {\rm Supp}(\phi))$ is a root of the 
polyhedralised discriminant $\tilde{\Delta}_{{\rm Supp}(\phi)}$.
\end{proposition}

\begin{proof}

Suppose that $P_0$ is a polyhedral root of $\phi$ of multiplicity at least two.  Since it is $\omega$-positive it has at least one vertex and so does $M \odot P_0$ (as usual, $M$ is the Minkowski sum of the coefficients of $\phi$) .
 Let $v$ be a vertex of $M \odot P_0$ and let $\ell$ be an element  in the interior of $N_{M \odot P_0}(v)$.  Consider the tropical polynomial $\phi_{\ell}(y):={\rm min}_{i \in {\rm Supp}(\phi)}(({\rm min}_{{\bf r} \in Q_i}  \ell({\bf r}))+i \cdot y)$.  Since minimising over $\ell$ commutes with respect to both $\oplus$ and $\odot$, the polynomial $\phi_{\ell}$ is also in product form.  Furthermore, $y_0={\rm min}_{{\bf r} \in P_0}(\ell({\bf r}))$ is a root of $\phi_{\ell}$ of multiplicity at least two. Hence, the tuple of coefficients of $\phi_{\ell}$ are a root of the  corresponding tropical discriminant \cite{DicFeiStu07}.  
 
 Note that for any $i \in {\rm Supp}(\phi)$, there exists a unique point that minimises $\ell$ over $Q_i$ for every $\ell \in N_{M \odot P_0}(v)$. Let this point be denoted by $q_i$.
 We claim that the tuple $(q_i \odot N^{\star}_{M \odot P_0}(v)|~i \in {\rm Supp}(\phi))$ is a root of $\tilde{\Delta}_{{\rm Supp}(\phi)}$.  
 

 
 To see this, consider the polyhedron obtained by evaluating  $\tilde{\Delta}_{{\rm Supp}(\phi)}$ at $(q_i \odot N^{\star}_{M \odot P_0}(v)|~i \in {\rm Supp}(\phi))$. Consider any vertex $w$ of this polyhedron. 
  
By construction, for every $\ell$ in the normal cone of $w$, the tuple $(\ell(q_i)|~i \in  {\rm Supp}(\phi))$ is a root of the tropical discriminant associated with ${\rm Supp}(\phi)$. Furthermore,  
there exists a full-dimensional subcone $C$ of the normal cone of $w$ and distinct elements $(\alpha_i|~i \in {\rm Supp}(\phi)),(\beta_i|~i \in {\rm Supp}(\phi))$ in the support of the tropical discriminant (associated with ${\rm Supp}(\phi)$)  such that for every $\ell \in C$ the corresponding minimum is attained at these elements.  More precisely,  $\sum_{i \in {\rm Supp}(\phi)} \alpha_i \cdot \ell(q_i)=\sum_{i \in {\rm Supp}(\phi)}  \beta_i \cdot \ell(q_i) \leq \sum_{i \in {\rm Supp}(\phi)} \gamma_i \cdot \ell(q_i)$ for all $(\gamma_i|~i \in {\rm Supp}(\phi))$ in the support of the tropical disriminant.  By Principle \ref{linalg_princ}, this implies that $w=\sum_{i \in {\rm Supp}(\phi)} \alpha_i \cdot q_i=\sum_{i \in {\rm Supp}(\phi)} \beta_i \cdot q_i$.  Hence, $w$ is shared by at least two summands of $\tilde{\Delta}_{{\rm Supp}(\phi)}(q_i \odot N^{\star}_{M \odot P_0}(v)|~i \in {\rm Supp}(\phi))$. We conclude that $(q_i \odot N^{\star}_{M \odot P_0}(v)|~i \in {\rm Supp}(\phi))$ is a root of the polyhedralised discriminant $\tilde{\Delta}_{{\rm Supp}(\phi)}$. 


Consider the evaluation $D:=\tilde{\Delta}_{{\rm Supp}(\phi)}(Q_i|~i \in {\rm Supp}(\phi))$ of the polyhedralised discriminant at $(Q_i|~i \in {\rm Supp}(\phi))$. Any vertex $w$ of $D$ is a vertex of  $E:=\tilde{\Delta}_{{\rm Supp}(\phi)}(q_i \odot N^{\star}_{M \odot P_0}(v)|~i \in {\rm Supp}(\phi))$ for suitable choices of $q_i$ and $v$. More precisely,  there is a full-dimensional subcone $C$ of $N_D(w)$ such that every $\ell \in C$ is uniquely minimised over each of the polyhedra $Q_i$ and $P_0$. The points $q_i$ and $v$ may be chosen to be this unique minimum over $Q_i$ and $P_0$ respectively. Furthermore, for any $k \in {\rm Supp}(\tilde{\Delta}_{{\rm Supp}(\phi)})$, the vertices of the $k$-th summands of $D$ and $E$ agree. Hence, $(Q_i|~i \in {\rm Supp}(\phi))$ is a root of $\tilde{\Delta}_{{\rm Supp}(\phi)}$. \end{proof}



For a cone in a Euclidean vector space, its normalised solid angle measure is the normalised volume of its intersection with the unit sphere. More generally, the normalised solid angle measure $\mu$ of a polyhedron is the normalised solid angle measure of its normal cone.

\begin{proposition}\label{discpoly_prop}
Let $\phi=\oplus_i (Q_i  \odot Y^i)$ be a polyhedral polynomial in product form. Suppose that the tuple of coefficients $(Q_i|~i \in {\rm Supp}(\phi))$ of $\phi$ is a root of $\tilde{\Delta}_{{\rm Supp}(\phi)}$, then there is a root of $\phi$ of multiplicity at least two whose normalised solid angle measure is at least $\delta \cdot \mu(M)$ where $\delta=\frac{1}{(\binom{|{\rm Supp}(\phi)|}{3} \cdot |V(M)|)}$.
\end{proposition}
\begin{proof}
Via an averaging argument, note that there is a vertex $v$ of $M$ such that $\mu(N_M(v)) \geq \mu(M)/|V(M)|$. Let $\sum_{i \in {\rm Supp}(\phi)}q_i$ be its Minkowski decomposition. For each $\ell \in N_M(v)$,  the tuple $(\ell(q_i)|~i \in {\rm Supp}(\phi))$ satisfies the tropical discriminant corresponding to ${\rm Supp}(\phi)$. Since $\phi$ is in product form, the tropical polynomial $\phi_\ell$ is also in product form for every $\ell \in N_M(v)$. Hence, each $\phi_\ell$ has a root $y_{\ell}$ such that the corresponding minimum is attained by at least three summands. Furthermore, there is a full-dimensional subcone $C$ of $N_M(v)$ and distinct elements $i_1,i_2,i_3 \in {\rm Supp}(\phi)$ such that $\mu(C) \geq \mu(N_M(v))/\binom{|{\rm Supp}(\phi)|}{3}$ and
\begin{equation}\label{triple_eq}
i_1 \cdot y_{\ell}+\ell(q_{i_1})=i_2 \cdot y_{\ell}+\ell(q_{i_2})=i_3 \cdot y_{\ell}+\ell(q_{i_3}) \leq i \cdot y_{\ell}+\ell(q_{i}) 
\end{equation}

for all $i \in {\rm Supp}(\phi)$.  Hence, $y_{\ell}=\ell(\rho_{i_1,i_2}^{(v)})=\ell(\rho_{i_2,i_3}^{(v)})=\ell(\rho_{i_1,i_3}^{(v)})$ for all $\ell \in C$. Since $C$ is full-dimensional, this implies that $\rho_{i_1,i_2}^{(v)}=\rho_{i_2,i_3}^{(v)}=\rho_{i_1,i_3}^{(v)}$. Furthermore, by Equation \ref{triple_eq}, the affine cone $\rho_{i_1,i_2}^{(v)} \odot C^{\star}$ is a root of $\phi$ of multiplicity at least two.  \end{proof}

At the time of writing, it is not clear if the ``local'' affine cone solutions in Proposition \ref{discpoly_prop} can be glued together to construct a root of higher multiplicity whose Minkowski-Weyl cone is equal to that of $M$. For this, the key additional condition that these local solutions need to satisfy is Condition \ref{comp2_con} in the definition of local compatible system.

\section{Appendix}

\subsection{Proof of Proposition \ref{prodform_prop}}

We start with a proposition on the recession cones of the roots of a polynomial. 

\begin{proposition} Let $p \in \mathbb{K}_P[x]$ be a non-zero polynomial with $\phi$ as its polyhedralisation. The convex set $\odot_{i}{\rm New}(\alpha_i)$ is a polyhedron and ${\rm MWC}(\odot_{i}{\rm New}(\alpha_i)) \subseteq  {\rm MWC}(M)$ where ${\rm MWC}(.)$ is the recession cone and $M$ is the Minkowski sum of the coefficients of $\phi$. 
 \end{proposition}

\begin{proof}
Write  $p/x^{k}=\prod_{\alpha_i \neq 0} (x-\alpha_i)$ where $k$ is the multiplicity of zero as a root of $p$ and $\alpha_i$ varies over its non-zero roots. Evaluating both sides at $0_{\mathbb{K}}$ and applying  ${\rm New}(.)$ yields $\odot_{i}{\rm New}(\alpha_i)={\rm New}(c_k)$ where $c_k$ is the coefficient of $x^k$ in $p$. Hence, $\odot_{i}{\rm New}(\alpha_i)$ is a polyhedron. The second part follows by applying {\rm MWC}(.) to the preceding equation and noting that ${\rm MWC}({\rm New}(c_k)) \subseteq {\rm MWC}(M)$.
\end{proof}

\begin{corollary}\label{poly_cor}
For any non-zero root $\alpha_i$ of $p$, the convex set ${\rm New}(\alpha_i) \odot {\rm MWC}(M)$ is a polyhedron.
\end{corollary}

{\bf  Hausdorff-Angle Metric:} The Hausdorff metric is a metric on the set of compact convex sets \cite{JohLin01} and hence, induces a metric on the set of convex polytopes. 
However, this metric does not readily extend to sets containing convex polyhedra such as (the set underlying) $\mathcal{A}$. 
 In order to define a metric on $\mathcal{A}$,  we start with a metric on the set of polyhedral cones. This metric depends on an inner product on the ambient vector space $W$.
Suppose that $\mathbb{S}^n$ is the unit sphere with respect to this inner product.  For polyhedral cones $C_1,C_2$, the metric $d_{\rm co}$ is defined as
$d_{\rm co}(C_1,C_2)=d_{\rm hau}(C_1 \cap \mathbb{S}^n,C_2 \cap \mathbb{S}^n)$ where $d_{\rm hau}$ is the Hausdorff metric on $W$. 

Using this, we define a metric on $\mathcal{A}$ as follows. For polyhedra $P_1,P_2$ in $\mathcal{A}$, the distance $d(P_1,P_2):=d_{\rm hau}(P_1^{\rm poly},P_2^{\rm poly})+d_{\rm co}(P_1^{\rm cone},P_2^{\rm cone})$ where $P_i^{\rm poly}$ and $P_i^{\rm cone}$ are the Minkowski-Weyl polytope and the Minkowski-Weyl cone of $P_i$, respectively. More concisely, this is the  product metric (in the $\ell_1$ sense) induced by the Hausdorff metric on the set of polytopes and the metric $d_{\rm co}$ on the set of polyhedral cones.  
The following propositions set up the properties needed to carry out the proof of Proposition \ref{prodform_prop}. 
\begin{proposition}\label{density_prop}
There is a dense subset $\mathcal{D}$  of $\mathcal{A}$ with respect to the topology induced by $d$ such that $\phi(P)=Q_d \odot (\odot_{i=1}^{d} (P \oplus P_i))$ for every $P \in \mathcal{D}$.
\end{proposition}
\begin{proof}
Consider an element $s \in \mathbb{K}$ such that ${\rm New}(s)=P$. If ${\rm New}(p(s)) \neq \phi(P)$, then $P$ must have a vertex of the form $\rho^{(v)}_{i,j}$ for some vertex $v$ of $M$ and $i,j \in {\rm Supp}(\phi)$. Hence,  the locus $\mathcal{D}_1=\{P \in \mathcal{A}| {\rm~for~every~} s \in \mathbb{K} {\rm~such~ that~} {\rm New}(s)=P,~{\rm New}(p(s))=\phi(P)\}$ is dense in $\mathcal{A}$. 
Similarly, if ${\rm New}(c_d \prod_{i=1}^{d} (s-\alpha_i)) \neq Q_d \odot (\odot_{i=1}^{d} (P \oplus P_i))$, then $P$ must share a vertex with ${\rm New}(\alpha_i)$ for some $i$ from $1$ to $d$.
Hence, the locus $\mathcal{D}_2=\{P \in \mathcal{A}|{\rm~for~every~} s \in \mathbb{K} {\rm~with~} {\rm New}(s)=P,~ {\rm New}(c_d \prod_{i=1}^{d} (s-\alpha_i))=Q_d \odot (\odot_{i=1}^{d} (P \oplus P_i)) \}$ is also dense in $\mathcal{A}$. The set $\mathcal{D}=\mathcal{D}_1 \cap \mathcal{D}_2$ and is hence, also dense in $\mathcal{A}$.
\end{proof}

\begin{proposition}\label{cont_prop}
Every polyhedral polynomial, i.e. an element of  $\mathcal{A}[Y]$, is a continuous function with respect to the topology induced by $d$. 
\end{proposition}

\begin{proof}
It suffices to show that if $f,g:\mathcal{A} \rightarrow \mathcal{A}$ are continuous functions, then so are $f \oplus g$ and $f \odot g$. This follows from the property of the metric that 
$d(P_1 \oplus P_3,P_2 \oplus P_3)$ and $d(P_1 \odot P_3 ,P_2 \odot P_3)$ are both upper bounded by $d(P_1,P_2)$ for all $P_1,P_2,P_3 \in \mathcal{A}$.
\end{proof}

The proof sketch of Proposition \ref{prodform_prop} can now be made rigorous by invoking Proposition \ref{density_prop} and Proposition \ref{cont_prop}. 

\subsection{Basic Properties of VCC}\label{vcc_app}

In the following, we provide proofs for the propositions on VCC stated without proof.

{\bf Proof of Proposition \ref{vccop_prop}:}  By the hypothesis on the supports, we deduce that the sets $\hat{G}$ and $\tilde{G}$ are non-empty.
 
In the following, we show that the cones $N_{\hat{G}}(v)$ and $N_{\tilde{G}}(v)$ are full-dimensional for all $v$. If $v \in G_1 \cap G_2$, then $N_{\hat{G}}(v)=N_{G_1}(v) \cap N_{G_2}(v)$ and is, by definition, full-dimensional.  Otherwise, if the vertex $v \in G_i \setminus G_j$, then by definition there is an element $\ell$  in the interior of $N_{G_i}(v)$ such that $\ell(v)<\ell(u)$ for all $u \in G_j$.  Hence, there is an open neighbourhood $B$ of $\ell$, that is contained in  $N_{G_i}(v)$, such that every element $\ell_a \in B$ satisfies $\ell_a(v)<\ell_a(u)$ for all $u \in G_j$.  Since $N_{\hat{G}}(v)=\{\ell \in N_{G_i}(v)|~\ell(v) \leq \ell(u)~\forall u \in G_j\}$,  we conclude that $B \subseteq N_{\hat{G}}(v)$ and hence, $N_{\hat{G}}(v)$ 
 is a full-dimensional cone.  The full-dimensionality of $N_{\tilde{G}}(v)$ is immediate from its definition. By construction, the normal cones $N_{\hat{G}}(v)$ and  $N_{\tilde{G}}(v)$  are polyhedral cones for all valid $v$.
  Next, we show that $\{(v,N_{\hat{G}}(v))\}_{v \in \hat{G}}$ and  $\{(v,N_{\tilde{G}}(v))\}_{v \in \tilde{G}}$ satisfy the inequality in the definition of vertex-cone collections. 

If $v \in G_1 \cap G_2$,   the normal cone of $N_{\hat{G}}(v)=N_{G_1}(v) \cap N_{G_2}(v)$  and hence, $\ell(v) \leq \ell(u)$ for any $\ell \in N_{\hat{G}}(v)$ and any $u \in G_1 \cup G_2$.  On the other hand, the normal cone $N_{\hat{G}}(v)$  of a vertex $v \in G_i \setminus G_j$ is the subcone of $N_{G_i}(v)$ of elements $\ell $ such that $\ell(v) \leq \ell(u)$ for all $u \in G_i \cup G_j$.  Hence, 
$\mathfrak{G}_1 \oplus \mathfrak{G}_2$ is a vertex-cone collection.

Consider the collection $\{(v,N_{\tilde{G}}(v))\}_{v \in \tilde{G}}$. Suppose that the Minkowski decomposition of $v \in \tilde{G}$ is $v_1+v_2$. Since $N_{\tilde{G}}(v)=N_{G_1}(v_1) \cap N_{G_2}(v_2)$, the inequality $\ell(v_1+v_2) \leq \ell(u_1+u_2)$ holds for all $\ell \in N_{\tilde{G}}(v)$ and for all $u_1 \in G_1,~u_2 \in G_2$. Hence,  $\mathfrak{G}_1 \odot \mathfrak{G}_2$  is a vertex-cone collection.  \qed

{\bf Proof of Lemma \ref{supp_lemma}:}   

 We first treat the case  $\mathfrak{G}_1 \oplus \mathfrak{G}_2$.  By definition, ${\rm Supp}(\mathfrak{G}_1 \oplus \mathfrak{G_2}) \subseteq {\rm Supp}(\mathfrak{G}_1)$.  Consider an element $\ell \in {\rm Int}(N_{G_1}(v)) \cap {\rm Int}(N_{G_2}(u))$ where $v \in G_1,u \in G_2$ (here ${\rm Int}(.)$ denotes the interior).   Suppose that $\ell$ attains a unique minimum over $G_1 \cup G_2$. Without loss of generality, we may assume that $\ell(v)<\ell(u)$. This implies that  $v \in \hat{G}$ and $\ell \in N_{\hat{G}}(v)$. Otherwise, if $\ell$ does not attain a unique minimum over $G_1 \cup G_2$, then $\ell(u)=\ell(v)$. If $u=v$, then $N_{G_1}(v) \cap N_{G_2}(u)$ is full-dimensional and is equal to $N_{\hat{G}}(v)$. If $u \neq v$, then any Euclidean open neighbourhood of $\ell$ (as an element in $W^{\star}$) contains points     $\tilde{\ell}_1$ and $\tilde{\ell}_2$ such that $\tilde{\ell}_1(v)<\tilde{\ell}_1(u)$ and  $\tilde{\ell}_2(u)<\tilde{\ell}_2(v)$. Hence, $u$ and $v$ are both in $\hat{G}$  and $\ell \in N_{\hat{G}}(v) \cap N_{\hat{G}}(u)$.  We conclude that $\cup_{v_1 \in G_1,v_2 \in G_2}({\rm Int}(N_{G_1}(v_1)) \cap {\rm Int}(N_{G_2}(v_2)))$ is contained in ${\rm Supp}(\mathfrak{G}_1 \oplus \mathfrak{G_2})$.  The set $\cup_{v_1 \in G_1,v_2 \in G_2}({\rm Int}(N_{G_1}(v_1)) \cap {\rm Int}(N_{G_2}(v_2)))$ is dense in ${\rm Supp}(\mathfrak{G}_1)$.  Furthermore, both ${\rm Supp}(\mathfrak{G}_1)$ and  ${\rm Supp}(\mathfrak{G}_1 \oplus \mathfrak{G_2})$  are closed in the Euclidean topology. Hence, we conclude that 
 ${\rm Supp}(\mathfrak{G}_1) \subseteq {\rm Supp}(\mathfrak{G}_1 \oplus \mathfrak{G_2})$.
 
 
 
 Consider the case $\mathfrak{G}_1 \odot \mathfrak{G}_2$. By definition,  ${\rm Supp}(\mathfrak{G}_1 \odot \mathfrak{G_2}) \subseteq {\rm Supp}(\mathfrak{G}_1) \cap  {\rm Supp}(\mathfrak{G}_2)$. On the other hand, suppose that $\ell \in {\rm Int}(N_{G_1}(v_1)) \cap {\rm Int}(N_{G_2}(v_2))$ for some $v_1 \in G_1, v_2 \in G_2$. This implies that ${\rm Int}(N_{G_1}(v_1)) \cap {\rm Int}(N_{G_2}(v_2))$ is full-dimensional and hence, the point $v_1+v_2 \in \tilde{G}$ and $\ell \in N_{\tilde{G}}(v_1+v_2)$.
 Hence, $\cup_{v_1 \in G_1,v_2 \in G_2}({\rm Int}(N_{G_1}(v_1)) \cap {\rm Int}(N_{G_2}(v_2)))$  is contained in ${\rm Supp}({\mathfrak{G}_1} \odot \mathfrak{G}_2)$. Note that $\cup_{v_1 \in G_1,v_2 \in G_2}({\rm Int}(N_{G_1}(v_1)) \cap {\rm Int}(N_{G_2}(v_2)))$ is dense in ${\rm Supp}(\mathfrak{G}_1) \cap  {\rm Supp}(\mathfrak{G}_2)$.  The sets ${\rm Supp}(\mathfrak{G}_1 \odot \mathfrak{G_2})$ and $ {\rm Supp}(\mathfrak{G}_1) \cap  {\rm Supp}(\mathfrak{G}_2)$ are both closed (in the Euclidean topology). Hence, we conclude that   ${\rm Supp}(\mathfrak{G}_1) \cap  {\rm Supp}(\mathfrak{G}_2) \subseteq {\rm Supp}(\mathfrak{G}_1 \odot \mathfrak{G_2})$.   This concludes the proof of the lemma. \qed

{\bf Proof of Proposition \ref{commassogp_prop}:}   The proofs of commutativity of both the operations and associativity of Minkowski sum are straight forward and hence, are left out.  
 We discuss the proof for associativity of the convex hull operation.  By Proposition \ref{vccop_prop}, $\mathfrak{G}_1 \oplus (\mathfrak{G}_2 \oplus \mathfrak{G}_3)$ and $(\mathfrak{G}_1 \oplus \mathfrak{G}_2) \oplus \mathfrak{G}_3$ are both vertex-cone collections.  By Lemma \ref{supp_lemma}, both these vertex-cone collections have the same support and  we claim that they are both equal to the vertex-cone collection $\mathfrak{G}_{1,2,3}$ defined as follows. Its vertex set $V(\mathfrak{G}_{1,2,3})$ is the union of the following three sets:
 \begin{enumerate}
\item $\{v \in V(\mathfrak{G}_1) \cap V(\mathfrak{G}_2) \cap V(\mathfrak{G}_3)|~\cap_{k=1}^{3} N_{\mathfrak{G}_k}(v) \text{~is full-dimensional}\}$. 
\item $\{v \in V(\mathfrak{G}_i) \setminus V(\mathfrak{G}_j \oplus \mathfrak{G}_k)|\\~\text{ $\exists$ a full-dimensional open ball $B$ in $N_{\mathfrak{G}_i}(v):\forall \ell \in B$, $\ell(v) \geq \ell(u)$ for all $u \in V(\mathfrak{G}_j \oplus \mathfrak{G}_k)$}\}$ for $i, j, k$ that are all distinct. 
\item $\{v \in V(\mathfrak{G}_i \oplus \mathfrak{G}_j) \setminus V(\mathfrak{G}_k)|\\~\text{ $\exists$ a full-dimensional open ball $B$ in $N_{\mathfrak{G}_i \oplus \mathfrak{G}_j}(v): \forall \ell \in B,$ $\ell(v) \geq \ell(u)$ for all $u \in V(\mathfrak{G}_k)$}\}$ for $i ,j,k$ that are all distinct.
 \end{enumerate}
 The corresponding normal cones are:
 \begin{itemize}
\item  In the first case: $N_{\mathfrak{G}_{1,2,3}}(v)=\cap_{k=1}^{3}N_{\mathfrak{G}_k}(v)$.

\item In the second and third cases: we may assume, without loss of generality, that $v \in V(\mathfrak{G}_i)$ and under this assumption $N_{\mathfrak{G}_{1,2,3}}(v)=\{ \ell \in N_{\mathfrak{G}_i}(v)|~\ell(v) \leq \ell(u) {\rm~holds} ~\forall u \in V(\mathfrak{G}_j) \cup V(\mathfrak{G}_k)\}$.
 \end{itemize} 
 The correctness of this description follows from the definition of convex hull and the hypothesis that the supports of $\mathfrak{G}_1,\mathfrak{G}_2$ and $\mathfrak{G}_3$ are the same. 
\qed

\footnotesize

\noindent {\bf Author's address:}

\smallskip
\noindent Department of Mathematics,\\
Indian Institute of Technology Bombay,\\
Powai, Mumbai, \\
India, 400076.
\end{document}